\documentclass[11pt,a4paper]{amsart}
\usepackage[T1]{fontenc}
\usepackage{lmodern}
\usepackage{amsmath,amssymb,amsthm,mathrsfs,bm}
\usepackage[a4paper,margin=30mm]{geometry}
\usepackage{enumitem}
\usepackage{microtype}
\usepackage[colorlinks=true,linkcolor=blue,citecolor=blue,urlcolor=blue]{hyperref}
\allowdisplaybreaks[2]
\numberwithin{equation}{section}
\newtheorem{theorem}{Theorem}[section]
\newtheorem{proposition}[theorem]{Proposition}
\newtheorem{lemma}[theorem]{Lemma}
\newtheorem{corollary}[theorem]{Corollary}
\theoremstyle{definition}

\theoremstyle{remark}
\newtheorem{remark}[theorem]{Remark}

\AtEndEnvironment{remark}{
	\unskip
	\hfill $\diamond$  
}

\makeatletter 
\@mparswitchfalse
\makeatother
\normalmarginpar 

\makeatletter
\patchcmd{\@settitle}
{\uppercasenonmath\@title}
{}{}{}
\patchcmd{\@settitle}
{\bfseries}
{\fontsize{16}{20}\selectfont\bfseries}
{}{}

\patchcmd{\@setauthors}
{\MakeUppercase{\authors}}
{\authors}{}{}
\patchcmd{\@setauthors}
{\footnotesize}
{\fontsize{11}{15}\selectfont}
{}{}

\renewcommand{\author@andify}{%
	\nxandlist
	{\unskip\quad}
	{\unskip\quad}
	{\unskip\quad}%
}

\patchcmd{\@setaddresses}
{\indent}
{\bigskip \noindent}
{}{}

\patchcmd{\@setaddresses}
{\indent\emailaddrname}
{\noindent\emailaddrname}
{}{}
\makeatother

\newcommand{\equ}{equation}
\newcommand{\R}{\mathbb{R}}
\newcommand{\C}{\mathbb{C}}
\def\mbs{\mathbb{S}}
\def\ig{\textit{g}}

\newcommand\ce{\mathcal{E}}
\newcommand\cf{\mathcal{F}}
\newcommand\cg{\mathcal{G}}
\newcommand\ch{\mathcal{H}}
\newcommand\ck{\mathcal{K}}
\newcommand\cl{\mathcal{L}}
\newcommand\cn{\mathcal{N}}
\newcommand\co{\mathcal{O}}
\newcommand\ct{\mathcal{T}}
\def\msc{\mathscr{C}}
\def\msd{\mathscr{D}}
\def\msh{\mathscr{H}}
\def\mss{\mathscr{S}}

\def\De{\Delta}
\def\ka{\kappa}
\def\al{\alpha}
\def\bt{\beta}
\def\de{\delta}
\def\Ga{\Gamma}
\def\la{\lambda}
\def\La{\Lambda}
\def\om{\omega}
\def\sa{\sigma}
\def\vr{\varepsilon}
\def\va{\varphi}
\def\pa {\partial}
\def\op{\oplus}
\def\ot{\otimes}
\def\id{\text{Id}}

\DeclareMathOperator{\tr}{tr}
\DeclareMathOperator{\spann}{span}
\DeclareMathOperator{\diver}{div}
\DeclareMathOperator{\sgn}{sgn}
\DeclareMathOperator{\vol}{vol}
\DeclareMathOperator{\dist}{dist}
\DeclareMathOperator{\supp}{supp}
\DeclareMathOperator{\diag}{diag}
\DeclareMathOperator{\spec}{Spec}
\DeclareMathOperator{\real}{Re}

\newcommand{\Crit}{\operatorname{Crit}}
\newcommand{\Conf}{\operatorname{Conf}}
\newcommand{\norm}[1]{\left\|#1\right\|}

\newcommand{\weakto}{\rightharpoonup}
\newcommand{\Hess}{\operatorname{Hess}}

\title[Energy gaps for critical Dirac equations]{Energy gaps for critical Dirac equations on spheres}
\author{Tian Xu}
\address{Tian Xu \newline Department of Mathematics, \newline Zhejiang Normal University, \newline Jinhua, Zhejiang, 321004, China}
\email{xutian@amss.ac.cn}

\begin{document}
\begin{abstract}
We study the critical Dirac equation
\[
 D\psi=H(x)|\psi|^{\frac{2}{m-1}}\psi
\]
on the round sphere $\mbs^m$, $m\ge2$, where $D$ is the Dirac operator.
For $H\equiv1$, the ground-state energy
$E_m=\frac{m^{m-1}}{2^{m+1}}\vol(\mbs^m)$ depends only on the dimension
$m$. In dimension two, the nonzero solution energies are exactly
$\pi,3\pi,4\pi,\ldots$. For $m\ge3$, we prove a gap
$(E_m,(1+\Ga_m)E_m)$ above the ground-state energy and give bounds for
its optimal relative size $\Ga_m$.
For positive $H$, let $\tau_H=E_m/(\max H)^{m-1}$ be the first bubble
energy. If $m\ge3$, every $H$ in a residual subset of the positive
$C^1$ functions admits solutions with infinitely many distinct weighted
energies in $(\tau_H,\min\{2,1+\Ga_m\}\tau_H)$ tending to $\tau_H$.
Along a subsequence these solutions form a single ground-state bubble
at a global maximum of $H$. For each $k\ge1$, smooth positive
coefficients admitting at least $k$ distinct energies in this interval
are dense among positive functions in $C^1$ and $W^{2,s}$,
$m/2<s<\infty$.
In the $C^2$ topology, however, an open set of positive coefficients
admits only the trivial solution. We approximate a smooth $H_*$ in
this set by smooth positive $H_j$ in every $C^{1,\al}$ and $W^{2,s}$
with $0<\al<1$ and $1\le s<\infty$. Solutions for $H_j$ form a single
bubble at the unique maximum of $H_*$ as their weighted energies tend
to $\tau_{H_*}$.
\end{abstract}
\maketitle

\section{Introduction}

Let $(\mbs^m,\ig_0)$ be the unit round sphere endowed
with its standard spin structure, where $m\ge2$. Write $\Sigma\mbs^m$
for its complex spinor bundle and $D=D_{\ig_0}$ for its Dirac operator.
For a spinor field $\psi$, let $|\psi|$ denote its pointwise Hermitian
norm. Given a positive function $H:\mbs^m\to\R$, we consider the equation
\begin{\equ}\label{eq:prescribed}
 D\psi=H(x)|\psi|^{\frac{2}{m-1}}\psi
 \qquad\text{on }\mbs^m.
\end{\equ}
The exponent
\[
 q=\frac{2m}{m-1}
\]
is the critical exponent for the embedding of $H^{1/2}(\mbs^m,\Sigma\mbs^m)$-spinors into
Lebesgue spaces. The conformal covariance of $D$ makes
\eqref{eq:prescribed} a natural spinorial analogue of the prescribed
scalar curvature equation. In dimension two, solutions without zeros
also have a geometric interpretation through the spinorial
Weierstra\ss\ representation, see for instance \cite{Friedrich1998,KusnerSchmitt1996,Taimanov1998}.

The analysis of \eqref{eq:prescribed} involves two familiar difficulties.
The associated functional is strongly indefinite, and the critical
Sobolev embedding is not compact. In particular, a bounded sequence of
solutions may concentrate. Isobe \cite{Isobe2011} established a global
compactness theorem which describes this loss of compactness by
Euclidean solutions of the limiting equation. The decay and regularity
of such solutions were investigated by Borrelli and Frank
\cite{BF2020}. In \cite{BMW2021}, Borrelli, Malchiodi and Wu classified
solutions at the ground-state energy as conformal images of Killing
spinors. And, more recently, Wang and Zhang \cite{WZ2026} obtained a
quantitative stability inequality for the sharp spinorial Sobolev
inequality.

Equation \eqref{eq:prescribed} leads us to consider the
constant-coefficient equation
\begin{\equ}\label{eq:standard}
	D\psi=|\psi|^{q-2}\psi
	\qquad\text{on }\mbs^m.
\end{\equ}
And a first natural question is ``Must a solution whose energy is sufficiently close to the ground-state
level already be a ground state?" The affirmative answer follows from
known results. For a sequence of solutions whose energies tend to the
ground-state level, global compactness \cite[Theorem~5.2]{Isobe2011}
leaves only one nonzero profile. After passing to a subsequence and
applying suitable conformal transformations, the sequence converges
strongly in $H^{1/2}(\mbs^m,\Sigma\mbs^m)$. The classification in
\cite[Theorem~1.2]{BMW2021} identifies the limit as a ground state, and
the nondegeneracy proved in \cite[Proposition~5.1]{Isobe2013} implies
that all solutions sufficiently close to this limit are themselves
ground states. A sharper question is then ``How far above the
ground-state energy does this rigidity persist?'' In dimension two, where
the ground-state energy is $E_2=\pi$, the geometric argument in the proof
of \cite[Theorem~1.3]{Xu2024} already suffices to exclude solution energies
in $(E_2,2E_2)$. Here we show that the
energies of nontrivial solutions are exactly
$\{E_2,\,3E_2,\,4E_2,\,5E_2,\dots\}$. The spinorial Weierstra\ss\
representation and the classification of finite-energy solutions of the
constant-mean-curvature system \cite{BC1985} associate to each such
solution a rational map of degree $d$, with the solution energy equal to
$dE_2$. Each ramification order is twice the vanishing order of the
spinor at the corresponding point. A degree-two map necessarily has a
ramification point of order one, which rules out the energy $2E_2$.
In higher dimensions, we use a lower bound for the gap to exclude
concentration in the variational construction for prescribed coefficients.

Set
\begin{\equ}\label{eq:constants}
 \om_m=\vol_{\ig_0}(\mbs^m),\qquad
 C_m=\left(\frac m2\right)^m\om_m,\qquad
 E_m=\frac{C_m}{2m}.
\end{\equ}
We write $(\cdot,\cdot)$ for the Hermitian product in the spinor
fibers and $\real(\cdot,\cdot)$ for its real part. We denote the
round volume form by $d\vol_{\ig_0}$.
The energy functional for \eqref{eq:standard} is
\begin{\equ}\label{eq:energy}
 \ce(\psi)=\frac12\int_{\mbs^m}(D\psi,\psi)\,d\vol_{\ig_0}
       -\frac1q\int_{\mbs^m}|\psi|^q\,d\vol_{\ig_0}.
\end{\equ}
Every nontrivial solution of \eqref{eq:standard} has energy at least $E_m$; equality defines a
ground-state solution. Our first result is as follows.

\begin{theorem}\label{thm:gap}
For every $m\ge2$, there exists a constant $\ka_m\in(0,1]$, depending
only on $m$, such that every nontrivial weak solution
$\psi\in H^{1/2}(\mbs^m,\Sigma\mbs^m)$ of \eqref{eq:standard}
satisfies 
\begin{\equ}\label{gap-alternative}
\text{either} \qquad \ce(\psi)=E_m
\qquad\text{or}\qquad
\ce(\psi)\ge(1+\ka_m)E_m.
\end{\equ}
In the first case, $\psi$ is a conformal image of a Killing spinor.

For $m=2$, the energies of nontrivial solutions are precisely
\[
 \big\{\pi d:d\in\{1,3,4,\ldots\}\big\}.
\]
In particular, the first excited energy is $3E_2=3\pi$.
\end{theorem}

In dimension two, the energy gap gives a gap in the critical values
of the conformally invariant spinorial quotient of
\cite[Section~2]{Ammann2009}
\[
 J(\va)=
 \frac{\displaystyle \Big(\int_{\mbs^2}|D\va|^{4/3}\,d\vol_{\ig_0}\Big)^{3/2}}
      {\displaystyle \int_{\mbs^2}(D\va,\va)\,d\vol_{\ig_0}},
\]
defined for $\va\in W^{1,4/3}(\mbs^2,\Sigma\mbs^2)$ with positive
denominator. The Euler-Lagrange equation of $J$ and the
invertibility of $D$ show that its critical points can be rescaled
to solve \eqref{eq:standard}. On nontrivial solutions of
\eqref{eq:standard}, $J(\psi)=2\sqrt{\ce(\psi)}$.
Thus $\bigl(2\sqrt{\pi},2\sqrt{3\pi}\bigr)$ contains no critical
value of $J$.

For a smooth metric $\ig$ on $\mbs^2$, let $\la_j^+(\ig)$ be its
$j$-th positive Dirac eigenvalue, counted with complex multiplicity,
and let $[\ig_0]$ be the smooth conformal class of $\ig_0$. Set
\[
 \La_j=\inf_{\ig\in[\ig_0]}
       \la_j^+(\ig)\vol_\ig(\mbs^2)^{1/2}.
\]
By conformal covariance, each nontrivial solution $\psi$ of \eqref{eq:standard} induces the generalized
metric $|\psi|^4\ig_0$, which has area $4\ce(\psi)$ and admits $1$ as a Dirac eigenvalue.
Rescaling this metric to unit area changes that eigenvalue to $2\sqrt{\ce(\psi)}=J(\psi)$.
The quaternionic structure of surface spinors gives
$\La_{2k-1}=\La_{2k}$ \cite[Section~1.1]{KMP2024}, and
\cite[Theorem~1.3]{Martynyuk2026} gives their common value
$2\sqrt{k\pi}$. On the disjoint union of $k$ unit spheres, the
first positive eigenvalue $1$ has complex multiplicity $2k$ and
yields the same normalized value. The exclusion of the solution
energy $2\pi$, together with a nodal estimate for weighted
eigenspinors, gives the nonattainment of $\La_3=\La_4$ even among
generalized conformal metrics (see Proposition~\ref{prop:second-pair-nonattainment} in Section \ref{subsec:sharp-two-dimensional-gap}).

\medskip

We next present an application of Theorem \ref{thm:gap} to the existence problem for 
\eqref{eq:prescribed}. Given $H\in C^2(\mbs^m)$, let
\[
 H_{\max}=\max_{\mbs^m}H,\qquad
 H_{\min}=\min_{\mbs^m}H,\qquad
 \ch=\{x\in\mbs^m:H(x)=H_{\max}\}.
\]
For $\de>0$, write
$\ch_\de=\{x:\dist_{\ig_0}(x,\ch)<\de\}$.
We write $\Crit(H)=\{x\in\mbs^m:\nabla H(x)=0\}$ for the critical set,
where $\nabla H$ is the gradient with respect to $\ig_0$.
We use the following condition from \cite{Xu2024}:
\begin{enumerate}[label=$(\mathrm{H}^*)$]
\item\label{hyp:topology}
The inclusion of $\ch$ into $\ch_\de$ is not null-homotopic for some
sufficiently small $\de>0$. There exists
\[
 d\in\bigl(\max\{2^{-1/(m-1)}H_{\max},H_{\min}\},H_{\max}\bigr)
\]
such that the inclusion of $\ch$ into $\{H\ge d\}$ is null-homotopic.
There is no critical value of $H$ in $(d,H_{\max})$. At every critical
point $x$ with $H(x)=d$, the Hessian of $H$ is positive definite.
\end{enumerate}
One can view \textup{($\mathrm{H}^*$)} as a pinching of the high
superlevel sets, where lowering the level to $d$ allows $\ch$ to
contract to a point within $\{H\ge d\}$. 
Our existence result reads as

\begin{theorem}\label{thm:existence}
Let $m\ge3$ and let $H\in C^2(\mbs^m)$ be positive. Suppose that
\textup{($\mathrm{H}^*$)} holds with a level $d$ satisfying
\begin{\equ}\label{eq:pinching}
 \left(\frac{H_{\max}}d\right)^{m-1}<1+\ka_m,
\end{\equ}
where $\ka_m$ is as in Theorem~\ref{thm:gap}. Then
\eqref{eq:prescribed} has a nontrivial weak solution. This solution
satisfies
\begin{\equ}\label{eq:existence-energy}
 \frac{C_m}{H_{\max}^{m-1}}
 <\int_{\mbs^m}H|\psi|^q\,d\vol_{\ig_0}
 <\frac{C_m}{d^{m-1}}.
\end{\equ}
\end{theorem}

The existence argument uses the part of the energy gap below
$2E_m$. Each nonzero component of the concentration decomposition
in \cite{Isobe2011} carries energy at least $E_m$, so below $2E_m$
at most one such component can occur. For \eqref{eq:prescribed},
the corresponding threshold is $2E_m/(\max H)^{m-1}$. We therefore
take $\ka_m$ in Theorem \ref{thm:gap} to be the largest constant in $(0,1]$ for which
\eqref{gap-alternative} holds, as specified in
\eqref{eq:kappa-canonical}. For $m\ge3$, the exact upper endpoint
remains open. The lower bound in Subsection~\ref{subsec:gap-bound}
suffices for the existence argument. It depends on the Sobolev
stability constant from \cite{WZ2026}, an elliptic estimate and a
local uniqueness radius, all of which can be chosen in terms of $m$.

Let $\msd_m$ be the class of positive smooth functions that are
constant or satisfy \textup{($\mathrm{H}^*$)} and \eqref{eq:pinching}.
For solutions of \eqref{eq:prescribed}, the weighted energy is given by
$(2m)^{-1}\int_{\mbs^m}H|\psi|^q\,d\vol_{\ig_0}$. We write
$\tau_H=E_m/H_{\max}^{m-1}$ for the energy of a ground-state bubble
with coefficient $H_{\max}$. To construct several solutions for the
same coefficient, we apply the local modification
\eqref{eq:volcano-modification} successively near the maximum of each
smooth Morse approximation. Each step introduces a new min-max
level closer to $\tau_H$, while
Lemma~\ref{lem:subcritical-minmax-persistence} keeps the preceding
levels in disjoint energy intervals.

\begin{corollary}
\label{cor:stronger-density}
Let $m\ge3$ and let $k\ge1$ be an integer. The coefficients
$H\in\msd_m$ for which
\eqref{eq:prescribed} has at least $k$ nontrivial solutions with
pairwise distinct weighted energies in
\[
 \bigl(\tau_H,(1+\ka_m)\tau_H\bigr)
\]
are dense in $C^1_+(\mbs^m)=\{H\in C^1(\mbs^m):\min H>0\}$ and in
$
 W^{2,s}_+(\mbs^m)=\{H\in W^{2,s}(\mbs^m):\min H>0\}$,
$m/2<s<\infty$.
Here the minimum is taken for the continuous representative.
Moreover, for every $0<\al<1$, each positive smooth coefficient can
be approximated in $C^{1,\al}$ by such coefficients.
\end{corollary}

The local modification can place the energy of a solution arbitrarily
close to $\tau_H$. By preserving the corresponding min-max level under
$C^1$-perturbations, we obtain the following conclusion for a fixed
coefficient in a residual set.

\begin{theorem}\label{thm:baire-multiplicity}
Let $m\ge3$. There is a residual subset of $C^1_+(\mbs^m)$, given by a
countable intersection of open dense sets in the $C^1$ topology, such
that for every $H$ in this subset, \eqref{eq:prescribed} has infinitely
many nontrivial weak solutions $\psi_j$ with pairwise distinct weighted
energies in $(\tau_H,(1+\ka_m)\tau_H)$ tending to $\tau_H$.
Moreover, $\norm{\psi_j}_\infty\to\infty$. After passing to a
subsequence, the solutions form a single ground-state bubble at a
global maximum of $H$.
\end{theorem}

\begin{remark}\label{rem:density-obstruction}
For $H\in C^2(\mbs^m)$ and a vector field $X$, write
$X(H)=\ig_0(\nabla H,X)$.
Every weak solution of \eqref{eq:prescribed} satisfies the spinorial
Kazdan-Warner identity \cite{AHA2007}
\begin{\equ}\label{eq:kw-obstruction}
 \int_{\mbs^m}X(H)|\psi|^q\,d\vol_{\ig_0}=0
 \qquad\text{for every conformal vector field }X.
\end{\equ}
Let $x_{m+1}$ denote the last ambient coordinate on $\mbs^m\hookrightarrow\R^{m+1}$.
For $H_0=c+a x_{m+1}$, $c>a>0$, the choice
$X=\nabla x_{m+1}$ rules out nontrivial solutions. This nonexistence
persists in a $C^2$ neighborhood of $H_0$; see
Proposition~\ref{prop:open-nonexistence}.
Thus the density statement cannot be strengthened to the $C^2$-topology, although approximation in the finite-exponent Sobolev norms
and in $C^{1,\al}$ for $\al<1$ is possible as stated above.
The local construction and its norm estimates are proved in
Subsection~\ref{subsec:stronger-density}. Our hypotheses give
constructive sufficient conditions, rather than a characterization of
all solvable coefficients.
\end{remark}

The variational construction used in Theorem \ref{thm:existence} is the
all-dimensional part of \cite{Xu2024}. In that paper, concentration (also called blow-up) is
excluded in dimension $2$ by a geometric rigidity argument. More
precisely, the limiting bubble is shown to have ground-state energy by
using the nodal estimate, the spinorial Weierstra\ss\ representation,
the Li-Yau inequality and Alexandrov's theorem. The role of
\eqref{eq:pinching} is to put the possible bubble in the smaller interval
excluded by Theorem~\ref{thm:gap}. This replaces that last geometric
step in dimension $2$. It does not require a classification of every solution with
energy below $2E_m$.

Theorem~\ref{thm:existence} gives a nontrivial spinor, with no assertion
about its zero set in higher dimensions. Its pinching condition concerns the high
superlevel region and imposes no global restriction on
$H_{\max}/H_{\min}$. Corollary~\ref{cor:large-contrast} gives admissible
coefficients with arbitrarily large values of this ratio.
The perturbative theory of \cite{Isobe2013,IX2023} only gives existence
and, in \cite[Theorem~1.1]{IX2023}, multiplicity for coefficients
close to a constant. Here, in Corollary~\ref{cor:stronger-density}, the target coefficient can
be any positive function in the stated spaces, and the solutions
are distinguished by their energies.

Recently, for $m\geq3$, the paper \cite{LT2026} studied blow-up for the critical Dirac equation in dimensions $m\ge3$, using the unweighted critical
mass. In our notation, the equation and the mass are
\[
D\psi_j=H_j|\psi_j|^{q-2}\psi_j \quad\text{on }\mbs^m \quad \text{and} \quad
A(\psi):=\int_{\mbs^m}|\psi|^q\,d\vol_{\ig_0},
\]
where $H_j$ is assumed to converge uniformly to a positive function $H$, as $j\to\infty$.
It is proved in \cite[Theorem~1.1]{LT2026} that, if $\limsup_j A(\psi_j)\le C_mH_{\max}^{-m}$ and if the solution sequence $\{\psi_j\}$ blows up, the entire limiting mass concentrates (up to a subsequence) to a
single ground-state bubble at a global maximum $a$ of $H$. In particular, the blow-up profile is isolated and simple if $H_j$ is $C^1$-uniformly bounded in a neighborhood of $a$ (here the terminology \textit{isolated} and \textit{simple} for a blow-up sequence are inherited from \cite{YYLi95} on the prescribing scalar curvature problem, see also \cite{KMS09} for related materials).  If, in addition, the coefficients $H_j$ converge
in $C^2$ around the point $a$, then the Hessian of $H$ must vanish at $a$. And consequently, when a function $H$ is given to satisfy $\diver_{\ig_0}\nabla H<0$ at every global maximum, for arbitrary $\widehat H$ in a $C^2$-neighborhood of this $H$, one can obtain a uniform $L^\infty$-bounds for the solution of the equation $D\psi=\widehat H|\psi|^{q-2}\psi$ on $\mbs^m$ with $A(\psi)<C_m\widehat H_{max}^{-m}+\vr$ for some $\vr>0$, see \cite[Theorem~1.3]{LT2026}.
The hypotheses of this compactness result can hold even when the
equation has only the trivial solution throughout a $C^2$-neighborhood.
Indeed, one may check $H_0(x)=1+\vr_0x_{m+1}$ with $0<\vr_0<1$ satisfies the
Laplacian condition. At its unique maximum
$p_N=(0,\ldots,0,1)$, we have $\Hess_{\ig_0}H_0(p_N)=-\vr_0\ig_0$
and $\diver_{\ig_0}\nabla H_0(p_N)=-m\vr_0<0$.
By our Proposition~\ref{prop:open-nonexistence}, the equation has only
the trivial solution for every coefficient in a sufficiently small
$C^2$-neighborhood of $H_0$.

When $H\equiv h>0$, the weighted energy is $hA(\psi)/(2m)$.
At a global maximum $a$ as above, the limiting bubble solves the
equation with constant coefficient $H(a)=H_{\max}$. Its unweighted
mass $C_mH_{\max}^{-m}$ therefore corresponds to the weighted energy
\[
 \frac{H_{\max}}{2m}\,C_mH_{\max}^{-m}
 =\frac{E_m}{H_{\max}^{m-1}}=\tau_H.
\]
This is the lower endpoint of the energy interval in
Theorem~\ref{thm:existence}, and it equals $E_m$ when $H\equiv1$.

Suppose that $0<d<H_{\max}$ and that $H$ has no critical values in
$(d,H_{\max})$. Corollary~\ref{cor:C1-window-compactness} in Section \ref{sec:existence} gives
compactness under $C^1$ convergence of $H_j$ to $H$ when the weighted
energies of the nontrivial solutions stay in a fixed closed subinterval of
\[
 \bigl(\tau_H,\min\{E_m/d^{m-1},(1+\ka_m)\tau_H\}\bigr).
\]
Indeed, any concentration in this range would take the form of a
single bubble carrying the entire limiting energy. The Kazdan-Warner
identity places it at a critical point of $H$. The upper bound forces
that point to be a global maximum, and the gap then forces the
limiting energy to equal $\tau_H$, below the prescribed interval.

Proposition~\ref{prop:weaker-convergence-blowup} shows what can happen
as the weighted energies approach $\tau_H$ from above. We construct
smooth positive $H_j$ converging to a smooth $H$ in every $C^{1,\al}$
with $0<\al<1$ and every $W^{2,s}$ with a finite $s$. The corresponding nontrivial solutions form a single ground-state bubble at the unique maximum
of $H$, with unweighted masses tending to $C_mH_{\max}^{-m}$. Their
blow-up is isolated and simple, as described in
\cite[Theorems~1.1 and~1.2]{LT2026}. The limit $H$ has strictly
negative Laplacian at its maximum, and the $C^2$ norms of $H_j$ stay
bounded. However, the Hessians of $H_j$ at that point stay a fixed
distance from the Hessian of $H$. This illustrates the role of $C^2$-closeness 
in the uniform estimate of \cite[Theorem~1.3]{LT2026}.

The conformal invariance of \eqref{eq:standard} and of the scalar
Yamabe equation on the round sphere in dimension $m\ge3$ leads to a
similar question about the energy above the ground state. In the
scalar case, the positive
and negative parts of a sign-changing solution each carry at least
one ground-state energy. Weth \cite{Weth2006} proved a uniform lower
bound strictly above twice the ground-state energy for all such
solutions. Spinors have no corresponding decomposition into positive and
negative parts. For $m\ge3$, whether every nontrivial solution of
\eqref{eq:standard} below $2E_m$ is a ground state remains open.
There is also a parallel for variable coefficients. For prescribed
scalar curvature, Li \cite[Corollary~0.16]{YYLi96} obtained density of
realizable coefficients, and Ndiaye \cite[Theorem~1.1]{Ndiaye2006}
showed that any positive $C^2$ coefficient can be approximated
in $C^{1,\al}$, $0<\al<1$, by coefficients giving arbitrarily many
positive solutions. In the spinorial setting,
\cite[Theorem~1.4]{ISX2024} constructed a fixed positive coefficient
close to a constant and a sequence of solutions whose $L^\infty$
norms tend to infinity.
Corollary~\ref{cor:stronger-density} gives a spinorial counterpart to
the density and multiplicity results by approximating any positive
coefficient in the stated topologies with one that gives at least as
many distinct solution energies as prescribed.

The paper is organized as follows. Section~\ref{sec:prelim} collects
notation and analytic facts. Section~\ref{sec:local} proves local
rigidity of the ground-state manifold and an explicit transverse
spectral estimate. Section~\ref{sec:gap} establishes the energy gap
and studies its size. Section~\ref{sec:existence} proves existence,
density and generic multiplicity in the $C^1$ topology, and gives a
concentration example.

\section{Notation, definitions and known results}\label{sec:prelim}

\subsection{Spinors and the variational framework}

We write $\nabla$ for the spin connection when it acts on spinors.
For tangent vectors $X,Y$ and a spinor $\psi$, our Clifford convention is
\[
 X\cdot Y\cdot\psi+Y\cdot X\cdot\psi=-2g_0(X,Y)\psi.
\]
Clifford multiplication by a tangent vector is skew-adjoint with
respect to the Hermitian product.
The Dirac operator has spectrum
\[
 \spec(D)=\left\{\pm\left(\frac m2+k\right):
 k=0,1,\ldots\right\}.
\]
In particular, it is invertible. Set
\[
 |D|=(D^2)^{1/2},\qquad \sgn D=D|D|^{-1},
\]
with these operators and their powers defined by spectral calculus.
We regard
\[
 \msh=H^{1/2}(\mbs^m,\Sigma\mbs^m)
\]
as a real Hilbert space with the equivalent norm
\[
 \norm{\psi}_{\msh}^2=
 \int_{\mbs^m}\bigl(|D|^{1/2}\psi,|D|^{1/2}\psi\bigr)\,d\vol_{\ig_0}.
\]
The decomposition into the positive and negative spectral subspaces is
written as $\msh=\msh^+\op\msh^-$.
For $\psi\in\msh$, its components in these subspaces are
denoted by $\psi^+$ and $\psi^-$. 
For $1\le s\le\infty$, we write $\norm{\cdot}_s$ for the $L^s$ norm on
$\mbs^m$ with respect to $d\vol_{\ig_0}$.

Throughout the paper,
\begin{\equ}\label{eq:dual-exponents}
 q=\frac{2m}{m-1},\qquad \hat q=\frac{2m}{m+1},\qquad
 \frac1q+\frac1{\hat q}=1.
\end{\equ}
The embedding $\msh\hookrightarrow L^q$ is continuous and is
compact for exponents strictly below $q$. Since $q>2$, the functional
\eqref{eq:energy} is of class $C^2$ on the real Hilbert space
$\msh$. Its critical points are the weak solutions of
\eqref{eq:standard}. All weak equations and derivatives below are
understood over the real numbers.

\begin{lemma}\label{lem:elliptic}
There exists a constant $c_m>0$ such that
\begin{\equ}\label{eq:elliptic}
 \norm{v}_{\msh}\le c_m\norm{Dv}_{\hat q}
\end{\equ}
whenever $Dv\in L^{\hat q}$. In particular, every $\msh$-weak solution of
\eqref{eq:standard} belongs to $W^{1,\hat q}$.
\end{lemma}
\begin{proof}
The embedding $\msh\hookrightarrow L^q$ gives, by duality,
$L^{\hat q}\hookrightarrow H^{-1/2}$. Since
$D:\msh\to H^{-1/2}$ is an isomorphism, this proves
\eqref{eq:elliptic}. If $\psi\in\msh$ is a weak solution, then
$|\psi|^{q-2}\psi\in L^{\hat q}$ because $(q-1)\hat q=q$. Global elliptic
regularity and the invertibility of $D$ yield $\psi\in W^{1,\hat q}$.
\end{proof}

Testing the equation with the solution gives
\begin{\equ}\label{eq:on-solutions}
 \int_{\mbs^m}(D\psi,\psi)\,d\vol_{\ig_0}
 =\int_{\mbs^m}|\psi|^q\,d\vol_{\ig_0},
 \qquad
 \ce(\psi)=\frac1{2m}\int_{\mbs^m}|\psi|^q\,d\vol_{\ig_0}.
\end{\equ}
For a constant $h>0$, the equation $D\psi=h|\psi|^{q-2}\psi$ has the energy functional
\[
 \ce_h(\psi)=\frac12\int_{\mbs^m}(D\psi,\psi)\,d\vol_{\ig_0}
       -\frac hq\int_{\mbs^m}|\psi|^q\,d\vol_{\ig_0}.
\]
The change $u=h^{(m-1)/2}\psi$ transforms this equation into
\eqref{eq:standard}, and
\begin{\equ}\label{eq:constant-scaling}
 \ce_h(\psi)=h^{-(m-1)}\ce(u).
\end{\equ}

For Euclidean space $\R^k$, write $\ig_{\R^k}$ for the standard metric
and $D_{\R^k}$ for its Dirac operator.
By a finite-energy solution on $\R^m$ we mean a distributional solution
in the homogeneous Sobolev space
$\dot H^{1/2}(\R^m,\Sigma\R^m)$, which embeds continuously into
$L^q(\R^m,\Sigma\R^m)$. The Euclidean energy $\ce_{\R^m}$ is
obtained from \eqref{eq:energy} by replacing the sphere and its volume
form by $\R^m$ and $d\vol_{\ig_{\R^m}}$. We use the same constants
$C_m$ and $E_m$ for this conformally equivalent problem.

\subsection{Conformal symmetry and ground states}

Write $\Conf^+(\mbs^m)$ for the group of orientation-preserving
conformal diffeomorphisms of $\mbs^m$. Let $F\in\Conf^+(\mbs^m)$,
with $F^*\ig_0=e^{2u}\ig_0$. Let $B_u$ be the canonical
fiberwise isometry from the spinor bundle of $F^*\ig_0$ to that of $\ig_0$.
After choosing a spin lift, define the conformal action by
\[
 \ct_F\psi=e^{(m-1)u/2}B_u(F^*\psi).
\]
The conformal covariance formula reads
\[
 D(\ct_F\psi)=e^{(m+1)u/2}B_u(F^*D\psi).
\]
Consequently,
\begin{\equ}\label{eq:conformal-invariants}
 \norm{\ct_F\psi}_q=\norm{\psi}_q,\qquad
 \norm{D(\ct_F\psi)}_{\hat q}=\norm{D\psi}_{\hat q},
 \qquad \ce(\ct_F\psi)=\ce(\psi).
\end{\equ}
We will not require it to preserve the particular positive definite norm chosen on $\msh$.

The Killing spinors with Killing constant $-1/2$ form the first
positive eigenspace of $D$. We denote by $\ck$ the set of spinors
of pointwise unit length in this eigenspace,
\[
\ck=\Big\{\xi\,:\,\nabla_X\xi=-\frac12X\cdot\xi
\text{ for every vector field }X,\quad |\xi|\equiv1\Big\}.
\]
The set $\ck$ is a finite-dimensional compact sphere.

Set $S_m=(m/2)\om_m^{1/m}$. The sharp spinorial Sobolev inequality on
the round sphere (see \cite{WZ2026}) states that
\begin{\equ}\label{eq:sharp-spinorial-sobolev}
 \norm{D\psi}_{\hat q}^2
 \ge S_m\int_{\mbs^m}(D\psi,\psi)\,d\vol_{\ig_0}
 \qquad\text{for every }\psi\in W^{1,\hat q}.
\end{\equ}
All nonzero spinors that attain equality in \eqref{eq:sharp-spinorial-sobolev}
form the cone
\begin{\equ}\label{eq:optimizer-cone}
 \co=\big\{\ct_F(a\xi):a>0,\ \xi\in\ck,\ F\in\Conf^+(\mbs^m)\big\}.
\end{\equ}
With $a_0=(m/2)^{(m-1)/2}$, the spinor $a_0\xi$ solves
\eqref{eq:standard} for every $\xi\in\ck$. Fixing this amplitude
in \eqref{eq:optimizer-cone} gives the family
\begin{\equ}\label{eq:ground-family}
 \cg=\big\{\ct_F(a_0\xi):\xi\in\ck,\ F\in\Conf^+(\mbs^m)\big\}.
\end{\equ}
The ground-state classification of \cite{BMW2021} identifies $\cg$
with all nontrivial solutions of \eqref{eq:standard} having energy $E_m$.
The explicit Euclidean bubble formulas show that $\cg$ is a smooth finite-dimensional
manifold; equivalently, one may use translations, dilations and the unit
spinor parameter after stereographic projection. This is the same
critical manifold used in \cite{Isobe2013,IX2023}.


We end this section with the following refinement of \eqref{eq:sharp-spinorial-sobolev},
which is due to \cite[Theorem~1.1]{WZ2026}. We will use it for the estimate in Subsection~\ref{subsec:gap-bound}.

\begin{proposition}\label{prop:stability}
There exists a constant $b_m>0$, depending only on $m$, such that 
\begin{\equ}\label{eq:stability}
 \norm{D\psi}_{\hat q}^2
 -S_m\int_{\mbs^m}(D\psi,\psi)\,d\vol_{\ig_0}
 \ge b_m\inf_{\phi\in\co}\norm{D(\psi-\phi)}_{\hat q}^2,
\end{\equ}
for every spinor $\psi\in W^{1,\hat q}$.
\end{proposition}

\section{Local rigidity of the ground-state manifold}\label{sec:local}

We use the nondegeneracy of the ground-state manifold for $\ce$ to
prove local uniqueness and give an explicit estimate for the Hessian
on the orthogonal complement of its kernel.
We use the nonnegative scalar Laplacian
$\De =-\diver_{\ig_0}\nabla$.

\subsection{The linearized kernel and local uniqueness}

\begin{lemma}[The linearized kernel]\label{lem:kernel}
Let $\xi\in\ck$ and $\Phi=a_0\xi$. For a variation $\zeta\in\msh$,
write $f_\zeta=\real(\xi,\zeta)$. The real linearization of
$D\psi-|\psi|^{q-2}\psi$ at $\Phi$ is
\begin{\equ}\label{eq:linearization}
 L\zeta=D\zeta-\frac m2\zeta-\frac m{m-1}f_\zeta\xi.
\end{\equ}
Its kernel consists of the spinors $\zeta$ satisfying
\begin{\equ}\label{eq:kernel}
 \zeta=\eta+f_\zeta\xi+\frac1{m-1}\nabla f_\zeta\cdot\xi,
 \quad \De f_\zeta=mf_\zeta,
 \quad D\eta=\frac m2\eta,
 \quad \real(\xi,\eta)=0.
\end{\equ}
In particular,
\begin{\equ}\label{eq:kernel-tangent}
 \ker L=T_\Phi\cg.
\end{\equ}
\end{lemma}
\begin{proof}
By elliptic regularity, each $\zeta\in\ker L$ is smooth. The
spherical calculation in \cite[(5.8) and (5.24)]{Isobe2013} gives
$\De f_\zeta=mf_\zeta$. For any smooth real-valued function $f$, the
Killing spinor equation gives
\begin{\equ}\label{eq:killing-products}
 D(f\xi)=\nabla f\cdot\xi+\frac m2f\xi,
 \qquad
 D(\nabla f\cdot\xi)=(\De f)\xi-\frac{m-2}{2}\nabla f\cdot\xi.
\end{\equ}
If $\De f=mf$, these identities show that
$P(f)=f\xi+(m-1)^{-1}\nabla f\cdot\xi$ satisfies $LP(f)=0$ and $\real\bigl(\xi,P(f)\bigr)=f$.
Taking $f=f_\zeta$, we obtain $\eta=\zeta-P(f_\zeta)\in\ker L$ and
$\real(\xi,\eta)=0$. In particular, using the expression of $L$ in \eqref{eq:linearization}, it follows that $D\eta=(m/2)\eta$.
Conversely, for any $\zeta$ of the form \eqref{eq:kernel}, both
$\eta$ and $P(f_\zeta)$ belong to $\ker L$. The identity
\eqref{eq:kernel-tangent} follows directly from \cite[Proposition~5.1]{Isobe2013}
by stereographic projection.
\end{proof}

\begin{remark}\label{rem:amplitude}
Fix $\xi\in\ck$ and consider the curve $a\mapsto a\xi$ in $\co$.
This curve varies the amplitude $a>0$ while keeping the unit Killing
spinor $\xi$ fixed, and it passes through $\Phi=a_0\xi$ at $a=a_0$.
Its tangent vector at $a=a_0$ is $\xi$, which we call the amplitude
direction. Since $L\xi=-\frac{m}{m-1}\xi$, this direction lies outside
$\ker L$. Thus the right-hand side of
\eqref{eq:kernel-tangent} is $T_\Phi\cg$, not $T_\Phi\co$.
Phase rotations, and quaternionic symmetries when present, are already
included in the variation of the unit Killing spinor.
\end{remark}

\begin{proposition}[Local uniqueness]\label{prop:local}
For every $\Phi=a_0\xi$, $\xi\in\ck$, there is a neighborhood
$U_\Phi$ in $\msh$ such that every critical point of $\ce$ in
$U_\Phi$ belongs to $\cg$.
\end{proposition}
\begin{proof}
Let $\nabla\ce$ denote the gradient with respect to the
$\msh$ scalar product, and let $A=\nabla^2\ce(\Phi)$ be its
Hessian at $\Phi$. Define the operator $K$ by
\[
 K\zeta=|D|^{-1}\left(\frac m2\zeta+\frac m{m-1}f_\zeta\xi\right).
\]
For the chosen inner product on $\msh$, the Riesz map from $\msh$
to its dual is $|D|$. Denote the duality pairing between $\msh^*$
and $\msh$ by $\langle\cdot,\cdot\rangle$. For all
$\zeta,\eta\in\msh$, we have
\[
\langle |D|A\zeta,\eta\rangle
=(A\zeta,\eta)_{\msh}
=\nabla^2\ce(\Phi)[\zeta,\eta]
=\langle L\zeta,\eta\rangle.
\]
Thus $|D|A=L$. Using $|D|^{-1}D=\sgn D$ and the definition of $K$,
we obtain
\[
A=|D|^{-1}L=\sgn D-K.
\]
The operator $K$ is compact: the zero-order coefficients in $L$ are
smooth and bounded, and the inclusion $\msh\hookrightarrow L^2$
is compact. It follows that $A$ is self-adjoint and Fredholm of index
zero. Lemma~\ref{lem:kernel} gives
$\ker A=T_\Phi\cg$.

Set $W=(T_\Phi\cg)^\perp$, with orthogonality in $\msh$.
Let $P_W:\msh\to W$ be the orthogonal projection.
Self-adjointness implies that $A$ maps $W$ to itself and that its
restriction to $W$ is invertible. Choose a local parametrization
$z(t)$ of $\cg$ with $z(0)=\Phi$. The map
$(t,w)\mapsto z(t)+w$, $w\in W$, gives local coordinates in
$\msh$. Consider
\[
 \cf(t,w)=P_W\nabla\ce(z(t)+w).
\]
Its derivative in $w$ at $(0,0)$ is $A|_W$, hence is invertible.
Since $\cf(t,0)=0$, the implicit function
theorem implies that, for each $t$ sufficiently close to $0$, the
normal equation $\cf(t,w)=0$ has a unique small solution, namely $w=0$.
Now let $\psi=z(t)+w$ be a critical point of $\ce$ near $\Phi=z(0)$.
Then $\cf(t,w)=P_W\nabla\ce(\psi)=0$, so $w=0$ and
$\psi=z(t)\in\cg$.
\end{proof}

\subsection{An explicit transverse estimate}\label{subsec:transverse-estimate}

The linear part of the local uniqueness argument admits an explicit
bound. This bound concerns the Hessian, rather than the energy gap.

\begin{proposition}\label{prop:normal-spectral-gap}
Let $\xi\in\ck$, $\Phi=a_0\xi$, and let
$A=|D|^{-1}L$ be the Hilbert Hessian at $\Phi$, with $L$ as in
\eqref{eq:linearization}. Then
\begin{\equ}\label{eq:explicit-transverse}
 \norm{A\zeta}_{\msh}\ge\frac1{2m+3}\norm{\zeta}_{\msh}
 \qquad\text{for }\zeta\perp_{\msh}\ker A.
\end{\equ}
The constant is uniform in $\xi$ (but is not asserted to be optimal).
\end{proposition}
\begin{proof}
Choose a real $L^2$-orthonormal basis of spherical harmonics. For a
basis function $f$ of degree $k\ge1$, we have
$\De f=\la_k f$, where $\la_k=k(m+k-1)$. The normalized spinors
\begin{\equ}\label{e1e2}
e_1=f\xi,\qquad e_2=\la_k^{-1/2}\nabla f\cdot\xi
\end{\equ}
are $L^2$-orthogonal, and we can denote $V_f=\spann_{\R}\{f\xi,\,\nabla f\cdot\xi\}$. For distinct basis functions $f$ and $g$ of positive degree,
\[
V_f
\perp_{L^2}
V_g.
\]
Indeed, we have that $\real(f\xi,g\xi)=fg$, $\real(\nabla f\cdot \xi,g\xi)=\real(f\xi,\nabla g\cdot \xi)=0$ and
\[
\int_{\mbs^m}\real(\nabla f\cdot\xi,\nabla g\cdot\xi)\,d\vol_{\ig_0}
=\int_{\mbs^m}f\De g\,d\vol_{\ig_0}=0.
\]

Given $e_1$ and $e_2$ as in \eqref{e1e2}, for $\zeta=c_1e_1+c_2e_2$ with $c_1,\,c_2\in\R$, we have
$f_\zeta=\real(\xi,\zeta)=c_1f$, so $f_\zeta\xi=c_1e_1$.
Set $a=m/2$ and $b=m/(m-1)$. By \eqref{eq:killing-products} and the expression of $L$ in
\eqref{eq:linearization}, the restricted operators $D\big|_{V_f}$ and $L\big|_{V_f}$ has the following representation with respect to the basis $\{e_1,\,e_2\}$
\begin{\equ}\label{eq:transverse-blocks}
	D_k=\begin{pmatrix}a&\sqrt{\la_k}\\
		\sqrt{\la_k}&1-a\end{pmatrix},
	\qquad
	L_k=\begin{pmatrix}-b&\sqrt{\la_k}\\
		\sqrt{\la_k}&-(m-1)\end{pmatrix}.
\end{\equ}
Since $\la_k=k(k+2a-1)$, the eigenvalues of $D_k$ are
$\frac12\pm\sqrt{(a-\frac12)^2+\la_k}$, namely $a+k$ and $-(a+k-1)$.

The constant harmonic gives $\R\xi\perp_{L^2}V_f$ for each $f$ of
positive degree. To proceed, let us define
\[
 \mss=\overline{\R\xi\oplus\big(\oplus_f V_f\big)}^{\,L^2},
\]
where the sum $\oplus_f$ runs over the chosen basis of spherical harmonics functions of positive degree.
Write $\mss^\perp$ for the real $L^2$-orthogonal complement of $\mss$.
The spaces $V_f$ and $\R\xi$ are invariant under $D$ and $L$. Since
they are spanned by Dirac eigenspinors, the spectral projections of
$D$ preserve $\mss$ and $\mss^\perp$. By spherical harmonic expansion, $u\xi\in\mss$ for every real-valued
$u\in L^2(\mbs^m)$. For $\zeta\in\mss^\perp$, take $u=f_\zeta$.
Then $\zeta\perp_{L^2}f_\zeta\xi$, which gives $\norm{f_\zeta}_2^2=0$. Hence, the last
term in \eqref{eq:linearization} vanishes, and we have that
\[
\left.L\right|_{\mss^\perp\cap\msh}
=D-\frac m2=D-a.
\]
Since the inner product on $\msh$ is induced by $|D|$, the spaces
$\R\xi$, $V_f$, and $\mss^\perp\cap\msh$ form an orthogonal
decomposition of $\msh$. Each of these spaces is invariant under $A$.

On $\R\xi$, since $D\xi=a\xi$ and $L\xi=-b\xi$, we have
$A\xi=-(b/a)\xi=-2\xi/(m-1)$. Clearly $2/(m-1)>1/(2m+3)$, hence \eqref{eq:explicit-transverse} holds on $\R\xi$. On $\mss^\perp\cap\msh$, we have
$A=|D|^{-1}(D-a)$, so a Dirac eigenspinor with eigenvalue $\la$
is an eigenvector of $A$ with eigenvalue $(\la-a)/|\la|$.
The eigenspinors with $\la=a$ lie in $\ker A$. For $\la=a+k$ with
$k\ge1$, we have $k/(a+k)\ge1/(a+1)$, while for $\la<0$ the absolute
value is at least $1$. Thus the nonzero eigenvalues of $A$ on
$\mss^\perp\cap\msh$ have absolute value at least $2/(m+2)$,
which justifies \eqref{eq:explicit-transverse} on this subspace.

For $\zeta\in V_f$, let $c=(c_1,c_2)\in\R^2$ be its coordinate vector with respect to
the basis $\{e_1,\,e_2\}$, that is, $\zeta=c_1e_1+c_2e_2$. In this basis the $\msh$ inner product has
matrix $|D_k|$, and $A|_{V_f}$ has matrix representation $|D_k|^{-1}L_k$. In particular, we have
\[
 \norm{\zeta}_{\msh}=\big\||D_k|^{1/2}c\big\|_{\R^2}
 \quad\text{and} \quad
 \norm{A\zeta}_{\msh}=\big\||D_k|^{-1/2}L_kc\big\|_{\R^2},
\]
where $\|\cdot\|_{\R^2}$ stands for the canonical Euclidean norm on $\R^2$.

For $k\ge2$, put $t_m=m-1+b$ and $d_k=\la_k-m=(k-1)(m+k)$. Then
$\tr L_k=-t_m$ and $\det L_k=-d_k<0$, so the least singular value
of $L_k$ is
\[
 \ell_k=\frac{\sqrt{t_m^2+4d_k}-t_m}{2}
       =\frac{2d_k}{t_m+\sqrt{t_m^2+4d_k}}
       \ge\frac{d_k}{t_m+\sqrt{d_k}}.
\]
Since $a+k$ is the largest eigenvalue of $|D_k|$, we have
\[
 \norm{A\zeta}_{\msh}
 \ge\frac{\norm{L_kc}_{\R^2}}{\sqrt{a+k}}
 \ge\frac{\ell_k}{\sqrt{a+k}}\norm{c}_{\R^2}
 \ge\frac{\ell_k}{a+k}\norm{\zeta}_{\msh}.
\]
Using $\sqrt{d_k}\le m+k$ and $b\le2$, we obtain
\begin{align*}
 \frac{\ell_k}{a+k}
 &\ge\frac{(k-1)(m+k)}{(a+k)(t_m+\sqrt{d_k})}
 \ge\frac{k-1}{t_m+m+k}\\
 &\ge\frac1{t_m+m+2}
 =\frac1{2m+1+b}\ge\frac1{2m+3}.
\end{align*}

For $k=1$, substituting $\la_1=m=b(m-1)$ into
\eqref{eq:transverse-blocks} gives $\det L_1=0$ and
$\tr L_1=-(b+m-1)$. Thus $-L_1$ is positive semidefinite of rank one.
In the $\{e_1,\,e_2\}$ coordinates, the vectors
\[
 \frac{1}{\sqrt{m+1}}\binom{\sqrt m}{1},
 \qquad
 \frac{1}{\sqrt{m+1}}\binom{1}{-\sqrt m}
\]
form an orthonormal eigenbasis of $D_1$ with eigenvalues $a+1$ and
$-a$, respectively. In this basis, $|D_1|^{-1}$ has matrix
$\diag((a+1)^{-1},a^{-1})$ and the diagonal entries of $-L_1$ are
nonnegative. Moreover, a standard computation shows that
\[
 -\tr(|D_1|^{-1}L_1)
 \ge\frac{-\tr L_1}{a+1}
 =\frac{b+m-1}{a+1}\ge1.
\]
The matrix $|D_1|^{-1}L_1$ has rank one, so its trace is its only
nonzero eigenvalue. Recall that $A=|D|^{-1}L$ is self-adjoint in $\msh$, we obtain
$\norm{A\zeta}_{\msh}\ge\norm{\zeta}_{\msh}$ for $\zeta\in V_f$ with
$\zeta\perp_{\msh}\ker A$. Now, combining all the estimates on the
$\msh$-orthogonal summands gives \eqref{eq:explicit-transverse}.
\end{proof}

\section{The energy gap and its size}\label{sec:gap}

\subsection{Isolation of the ground-state energy}

\begin{proposition}
\label{prop:compactness-mod-conf}
Let $\psi_j\ne0$ be weak solutions of \eqref{eq:standard} satisfying
\[
 \sup_j\ce(\psi_j)\le\La<2E_m.
\]
After passing to a subsequence, there are orientation-preserving
conformal actions $\ct_j$ and a nontrivial solution $\psi_*$ of
\eqref{eq:standard} such that  $\ct_j\psi_j\to\psi_*$  in $\msh$.
In particular, the set of nontrivial solutions in every such energy
sublevel is sequentially compact modulo conformal transformations.
\end{proposition}
\begin{proof}
After taking a subsequence, write
$\ce(\psi_j)\to c\in[E_m,\La]$. The lower energy bound for
nontrivial solutions, on both the sphere and Euclidean space, is the
classical sharp bound recalled in \cite[Section 4]{Isobe2011}.
We first give an estimate specific to nontrivial solutions. For every
conformal action $\ct$, the spinor $\ct\psi_j$ again solves
\eqref{eq:standard}; hence Lemma~\ref{lem:elliptic} and conformal
invariance give
\begin{\equ}\label{eq:exact-conformal-H-bound}
 \norm{\ct\psi_j}_{\msh}
 \le c_m\norm{D(\ct\psi_j)}_{\hat q}
 =c_m\left(\int_{\mbs^m}|\psi_j|^q\,d\vol_{\ig_0}\right)^{1/\hat q}
 \le c_m(2m\La)^{1/\hat q}.
\end{\equ}
This is a uniform $\msh$-bound for the transformed solutions. 

Clearly,  $\{\psi_j\}$ is  a bounded Palais-Smale sequence for the functional $\ce$.
The global compactness theorem of \cite{Isobe2011} gives a weak
solution $\psi_\infty$ and finitely many nontrivial finite-energy
Euclidean solutions $\Psi^\ell$ such that
\begin{\equ}\label{eq:mod-conf-energy-decomposition}
 c=\ce(\psi_\infty)
     +\sum_{\ell=1}^k\ce_{\R^m}(\Psi^\ell).
\end{\equ}
Every nonzero term on the right-hand side is at least $E_m$. Since $c<2E_m$, exactly one
nonzero term occurs. If $\psi_\infty\ne0$, there is no blow-up, and
global compactness already gives strong convergence in $\msh$.
We then take $\ct_j$ to be the identity.

If $\psi_\infty=0$, the sequence $\{\psi_j\}$ develops a single bubble. We can choose
the conformal transformations using the center and scale of this
bubble. In stereographic coordinates, we translate the concentration center
to the origin and rescale so that the concentration scale is one.
These translations and dilations extend to conformal diffeomorphisms
of $\mbs^m$. Their actions $\ct_j$ on spinors give precisely the
rescaling used to extract the bubble, including the conformal
weight. Hence $\ct_j\psi_j$ has a nonzero weak limit on a fixed
coordinate ball.

By \eqref{eq:exact-conformal-H-bound}, we may pass to a subsequence
whose weak limit is denoted by $\psi_*$ in $\msh$. Applying the global compactness \cite{Isobe2011} to
$\ct_j\psi_j$, we obtain the strong convergence to $\psi_*$ in $\msh$.
Indeed, $\ce(\psi_*)\ge E_m$, and any further bubble would bring
the limiting energy to at least $2E_m$, contradicting $c<2E_m$.
In both cases, strong convergence gives $\ce(\psi_*)=c$.
\end{proof}

\begin{remark}\label{rem:exactness-mod-conf}
The argument uses nontrivial solutions in
\eqref{eq:exact-conformal-H-bound} and in applying global compactness
to the normalized sequence. It is not a statement about arbitrary
Palais-Smale sequences under unbounded conformal actions. The level
$2E_m$ enters through energy counting. We also point it out that no assertion is made that this
threshold is sharp for general nontrivial solutions.
\end{remark}

\begin{lemma}\label{lem:normalization}
Suppose that $\psi_j$ are nontrivial weak solutions of
\eqref{eq:standard} and $\ce(\psi_j)\to E_m$. After passing to a
subsequence, there exist conformal actions $\ct_j$ such that
$\ct_j\psi_j\to a_0\xi$ in $\msh$ for some $\xi\in\ck$.
\end{lemma}
\begin{proof}
For all large $j$, the energies are bounded by a fixed
$\La<2E_m$. Proposition~\ref{prop:compactness-mod-conf} gives
strong convergence, modulo conformal transformations, to a solution
of energy $E_m$. By the ground-state classification
\cite[Theorem~1.2]{BMW2021}, the limit belongs to $\cg$.
Let $\widetilde{\ct}_j$ be the conformal actions for which
$\widetilde{\ct}_j\psi_j\to\psi_*$ strongly in $\msh$.
Write $\psi_*=\ct_F(a_0\xi)$ and set
$\ct_j=\ct_F^{-1}\circ\widetilde{\ct}_j$.
Since the fixed inverse action $\ct_F^{-1}$ is bounded on $\msh$, one has
$\ct_j\psi_j\to a_0\xi$ in $\msh$.
\end{proof}

\begin{proof}[Proof of Theorem~\ref{thm:gap}]
If the conclusion were false, there would be a sequence of non-ground-state nontrivial
solutions $\psi_j$ with $\ce(\psi_j)\to E_m$. By
Lemma~\ref{lem:normalization}, their conformal transforms converge
strongly in $\msh$ to $a_0\xi$ for some $\xi\in\ck$.
Local uniqueness, Proposition~\ref{prop:local}, puts all sufficiently
large transformed solutions in $\cg$, contradicting their energies.
Thus the ground-state energy is isolated, and a positive constant
$\ka_m\le1$ with the asserted alternative exists.
\end{proof}

\begin{corollary}\label{cor:euclidean-gap}
Every nontrivial finite-energy solution of
\[
 D_{\R^m}\Psi=|\Psi|^{q-2}\Psi\quad\text{on }\R^m
\]
satisfies either
\[
 \int_{\R^m}|\Psi|^q\,d\vol_{\ig_{\R^m}}=C_m
 \quad\text{or}\quad
 \int_{\R^m}|\Psi|^q\,d\vol_{\ig_{\R^m}}\ge(1+\ka_m)C_m.
\]
\end{corollary}
\begin{proof}
By stereographic projection, a finite-energy Euclidean solution
corresponds to a solution on the punctured sphere. The removable
singularity result for finite-energy critical Dirac solutions extends
it across the missing point; this is part of the standard sphere-Euclidean
correspondence in \cite{BF2020,BMW2021}. The proof is completed by using the conformally invariance of the critical mass.
\end{proof}

To describe the size of the gap above $E_m$, define the set of
excited solutions by
\begin{\equ}\label{eq:least-excited-definition1}
 \msc_{\mathrm{exc}}(m)
 =\big\{\psi\in\msh:\nabla \ce(\psi)=0,\ \ce(\psi)>E_m\big\},
\end{\equ}
We denote the infimum of their energies by $E_{\mathrm{exc}}(m)$
and the optimal relative gap by $\Ga_m$, so
\begin{\equ}\label{eq:least-excited-definition2}
E_{\mathrm{exc}}(m)
=\inf_{\psi\in\msc_{\mathrm{exc}}(m)}\ce(\psi)\quad \text{and} \quad
\Ga_m=\frac{E_{\mathrm{exc}}(m)}{E_m}-1.
\end{\equ}
Here we adopt the convention that the infimum over the empty set is
$+\infty$. Theorem~\ref{thm:gap} gives $\Ga_m>0$. By definition,
$(E_m,(1+\Ga_m)E_m)$ is the largest open interval with left endpoint
$E_m$ containing no solution energy. The constant used in
Theorem~\ref{thm:existence} will be
\begin{\equ}\label{eq:kappa-canonical}
 \ka_m=\min\{1,\Ga_m\},
\end{\equ}
which is the truncated gap we settle for here, given the difficulty
of determining $\Ga_m$ exactly for $m\ge3$. We nevertheless obtain
a lower bound for $\Ga_m$ in terms of analytic constants in
Proposition~\ref{prop:safe-gap} and explicit upper bounds in
Proposition~\ref{prop:explicit-gap-upper}. Together these estimates
give a rough indication of the size of $\Ga_m$.

\subsection{The optimal gap in dimension two}
\label{subsec:sharp-two-dimensional-gap}

In dimension two, the spinorial Weierstra\ss\ representation relates
solutions of the Dirac equation to branched conformal immersions of
constant mean curvature. This correspondence gives the full set of
nonzero solution energies. Here $C_2=4\pi$ and $E_2=\pi$.

\begin{proposition}[Energy levels on $\mbs^2$]
\label{prop:two-dimensional-gap}
The energies of the nontrivial weak solutions
$\psi\in H^{1/2}(\mbs^2,\Sigma\mbs^2)$ of
\begin{\equ}\label{eq:two-dimensional-critical}
 D\psi=|\psi|^2\psi
\end{\equ}
form the set
\begin{\equ}\label{eq:two-dimensional-energy-quantization}
 \big\{\pi d\,:\,d\in\{1,3,4,5,\ldots\}\big\}.
\end{\equ}
Consequently,
\begin{\equ}\label{eq:optimal-two-dimensional-gap}
 E_{\mathrm{exc}}(2)=3E_2=3\pi,
 \qquad \Ga_2=2,\qquad \ka_2=1.
\end{\equ}
In particular, every nontrivial solution with $\ce(\psi)<3E_2$ belongs
to $\cg$.
\end{proposition}

\begin{proof}
Let $\psi\in H^{1/2}(\mbs^2,\Sigma\mbs^2)$ be a nontrivial weak solution of \eqref{eq:two-dimensional-critical}. Our proof is divided into the following four steps.

\smallskip\noindent\textit{Step 1: A globally defined branched immersion.}
By the Sobolev embedding, $\psi\in L^4$. The regularity theorem for
weak solutions of the cubic Dirac equation \cite{Wang2010} then
gives $\psi\in C^\infty$.
Since $\psi$ satisfies $D\psi=P\psi$ with the smooth potential
$P=|\psi|^2$, B\"ar's nodal-set theorem gives a discrete, hence
finite, zero set $Z$ on $\mbs^2$. The strong unique continuation
argument in its proof gives finite vanishing order at each zero, see
\cite{Baer1997}.

The spinorial Weierstra\ss\ representation
\cite[Section~9]{Ammann2009} constructs a closed,
$\R^3$-valued smooth one-form $\al_\psi$, quadratic in $\psi$.
Since $\mbs^2$ is simply connected, the closed one-form $\al_\psi$
is exact. Thus there is a smooth map $F:\mbs^2\to\R^3$, unique up
to translation, such that $dF=\al_\psi$. On $\mbs^2\setminus Z$
this map is a conformal immersion with mean curvature one, in the
convention where mean curvature is the average of the principal
curvatures. Moreover,
\begin{\equ}\label{eq:weierstrass-induced-metric-two-d}
 F^*\ig_{\R^3}=|\psi|^4\ig_0.
\end{\equ}
The branching index is the vanishing order of $dF$. Because
$\al_\psi$ is quadratic and has length proportional to $|\psi|^2$,
a zero of $\psi$ of order $k$ gives a branch of $F$ of order $2k$.
Thus every branching index is even.

\smallskip\noindent\textit{Step 2: The mean curvature equation.}
Choose an arbitrary point $p_\infty\in\mbs^2$ and a stereographic projection
$\varpi:\mbs^2\setminus\{p_\infty\}\to\R^2$, and set
$\om=F\circ\varpi^{-1}:\R^2\to\R^3$. On the regular set of $\om$, let $H_\om$
be the signed mean curvature with respect to the unit normal in
the direction of $\om_x\wedge\om_y$, where $\wedge$ is the cross
product in $\R^3$. In conformal coordinates, the ordinary coordinate
Laplacian of $\om$ equals twice its mean curvature vector times the
area density $|\om_x\wedge\om_y|$. This gives the equation
$\pa_x^2\om+\pa_y^2\om=2H_\om\,\om_x\wedge\om_y$, known as the
$H$-system \cite[p.~149, equation~(H)]{BC1984}. After replacing
$\om$ by $-\om$ if necessary, we have $H_\om=1$ and
\begin{\equ}\label{eq:two-dimensional-hsystem}
 \pa_x^2\om+\pa_y^2\om=2\om_x\wedge\om_y.
\end{\equ}
Both sides are smooth, so the equation holds across the branch
points as well.

Write $|dF|_{\mathrm{HS}}$ for the Hilbert-Schmidt norm of $dF$.
Conformal invariance of the Dirichlet integral and
\eqref{eq:weierstrass-induced-metric-two-d} give
\begin{\equ}\label{eq:spinor-hsystem-energy}
 \int_{\R^2}|\nabla\om|^2\,d\vol_{\ig_{\R^2}}
 =\int_{\mbs^2}|dF|_{\mathrm{HS}}^2\,d\vol_{\ig_0}
 =2\int_{\mbs^2}|\psi|^4\,d\vol_{\ig_0}<\infty.
\end{\equ}
Thus $\om\in\dot H^1(\R^2;\R^3)$, and the classification of Brezis-Coron
\cite[Lemma~0.1 and Appendix, Lemma~A.1]{BC1985}, also stated as
\cite[Theorem~2]{GLZ2025}, applies to
\eqref{eq:two-dimensional-hsystem}. Identify $\R^2$ with $\C$ and write
$\widehat{\C}=\C\cup\{\infty\}$ for the Riemann sphere. Let $\varpi_0$
be the standard stereographic projection of the unit sphere
$\mbs^2\subset\R^3$, extended by sending its North Pole to $\infty$. We regard its inverse as a map
$\varpi_0^{-1}:\widehat{\C}\to\R^3$. The classification gives a rational
map $R:\widehat{\C}\to\widehat{\C}$ of degree $d\ge1$ and a translation
vector $b\in\R^3$ such that
\[
 \om(z)=\varpi_0^{-1}(R(z))+b,
 \qquad
 \int_{\R^2}|\nabla\om|^2\,d\vol_{\ig_{\R^2}}=8\pi d.
\]
The degree is positive because $\psi\ne0$, so $F$ is not constant. Using
\eqref{eq:on-solutions} and \eqref{eq:spinor-hsystem-energy}, we obtain
\[
 \ce(\psi)=\frac14\int_{\mbs^2}|\psi|^4\,d\vol_{\ig_0}
 =\frac18\int_{\R^2}|\nabla\om|^2\,d\vol_{\ig_{\R^2}}=\pi d.
\]

\smallskip\noindent\textit{Step 3: Exclusion of degree two.}
Suppose $d=2$. The Riemann-Hurwitz formula gives total ramification
$2d-2=2$, so $R$ has a ramification point. At any such point its local
degree is at most $2$ and hence equals $2$. Its ramification order is
therefore $1$. Since $\varpi_0^{-1}$ is an immersion, the corresponding branch point of $F$ also has order $1$, contradicting the evenness of the branching indices established in Step~1. Thus $d\ne2$.
When $d=1$, the energy is $E_2$, and ground-state rigidity yields
$\psi\in\cg$. Every other solution consequently has energy at least
$3E_2$.

\smallskip\noindent\textit{Step 4: Realization of the energy levels.}
For odd $d$, take $R_d(z)=z^d$, as in \cite[Remark, p.~449]{Ammann2009}.
Since $R_d'(z)=dz^{d-1}$, the ramification order at $0$ is $d-1$,
and every other finite point is unramified. In the reciprocal
coordinate $w=1/z$, we have $1/R_d(1/w)=w^d$, so the ramification
order at $\infty$ is also $d-1$. Both orders are even.

For even $d\ge4$, set $t=d/(d-2)$. Taking $k=1$ in
\cite[Proposition~3]{ABEGKM2018} gives the normalized Belyi map
$f_d(z)=z^{d-1}(z-t)/(1-tz)$. Rescaling the source and the target,
we use
\[
 R_d(z)=-t^d f_d(z/t)=z^{d-1}\frac{z-t^2}{z-1}.
\]
A direct calculation shows that
\[
 R_d'(z)=(d-1)z^{d-2}\frac{(z-t)^2}{(z-1)^2}.
\]
Since $t>1$, the numerator and denominator have no common zero,
so $R_d$ has degree $d$. Away from the poles, the derivative gives
exactly two ramification points, $0$ and $t$, with orders $d-2$ and
$2$, respectively. In the reciprocal target coordinate, $1/R_d$
has a simple zero at $1$ and a zero of order $d-1$ at $\infty$.
The ramification orders at these two points are therefore $0$
and $d-2$, respectively. The orders sum to
$(d-2)+2+0+(d-2)=2d-2$, in agreement with the Riemann-Hurwitz
formula. Since $d$ is even, the ramification order at each point
is even.

Regard $R_d$ as a map between round spheres and compose it with
the inclusion of the target unit sphere into $\R^3$. As in Step~2,
we may replace the resulting branched immersion by its negative
so that its signed mean curvature is one. Since all branching
indices are even, the converse spinorial Weierstra\ss\
representation \cite[Section~9, p.~464]{Ammann2009} gives a smooth
solution $\psi_d$ of \eqref{eq:two-dimensional-critical} in the
unique spin structure on $\mbs^2$. The induced metric is
$|\psi_d|^4\ig_0=R_d^*\ig_0$. Since $R_d$ covers the unit sphere $d$
times, we obtain
\[
 \ce(\psi_d)
 =\frac14\int_{\mbs^2}|\psi_d|^4\,d\vol_{\ig_0}
 =\frac14\int_{\mbs^2}R_d^*d\vol_{\ig_0}
 =\pi d.
\]
The first energy above $E_2$ is therefore $3\pi$, which gives
\eqref{eq:optimal-two-dimensional-gap}.
\end{proof}

\begin{remark}
The truncation in \eqref{eq:kappa-canonical} was introduced for the
existence argument in Theorem~\ref{thm:existence}. The two-dimensional
result shows that the retained interval can be strictly shorter
than the optimal gap, since $\ka_2=1$ and $\Ga_2=2$.
There is no nontrivial solution at energy $2E_2$. This does not rule
out a Palais-Smale sequence consisting of two separating
ground-state bubbles with limiting energy $2E_2$.

In Step~1, a zero of $\psi$ of order $k$ gives a branch point of
$F$ with branching index $2k$. To apply the converse spinorial
Weierstra\ss\ representation in Step~4, we choose $R_d$ with even
ramification order at every branch point.
For even $d$, we start with
$R_a(z)=z^{d-1}(z-a)/(z-1)$, where $a\in\C\setminus\{0,1\}$.
The ramification orders at $0$ and $\infty$ are already $d-2$.
The remaining critical points of $R_a$ are the roots of
$(d-1)z^2-((d-2)a+d)z+(d-1)a$.
To give them even ramification order, we have to require the two roots to
coincide. The discriminant is
\[
 \big((d-2)a+d\big)^2-4(d-1)^2a
 =(a-1)\big((d-2)^2a-d^2\big).
\]
The choice $a=1$ cancels a factor and lowers the degree. For $d>2$,
the other choice is $a=(d/(d-2))^2=t^2$, as in Step~4.
When $d=2$, the discriminant reduces to $-4(a-1)$, so the only candidate is $a=1$. Putting this into $R_a$ and canceling the common factor $z-1$ then gives the degree-one map $R_a(z)=z$. This agrees with Step~3, where the Riemann-Hurwitz formula rules out degree-two maps with even ramification at every point. For even $d\ge4$, the only parameter in this family that preserves degree $d$ and gives even ramification at every point is $a=(d/(d-2))^2=t^2$.
\end{remark}

\begin{remark}[Explicit solutions on $\R^2$]
In stereographic coordinates the maps $R_d$ in Step~4 also give
explicit solutions on $\R^2$. Write $R_d=P_d/Q_d$ with the coprime
polynomials used above. The computation in Step~4 gives
$P_d'Q_d-P_dQ_d'=A_d^2$, where
\[
 A_d(z)=
 \begin{cases}
 \sqrt d\,z^{(d-1)/2}, & d\text{ odd},\\
 \sqrt{d-1}\,z^{(d-2)/2}(z-t), & d\ge4\text{ even}.
 \end{cases}
\]
In the flat trivialization in which
$D_{\R^2}(u,v)=(-2\pa_zv,\,2\pa_{\bar z}u)$, set
\[
 \Psi_d(z)=\frac{\sqrt{2}}{|P_d(z)|^2+|Q_d(z)|^2}
 \begin{pmatrix}
 A_d(z)\overline{Q_d(z)}\\
 -P_d(z)\overline{A_d(z)}
 \end{pmatrix}.
\]
The denominator is positive everywhere, so $\Psi_d$ is smooth even
at the poles of $R_d$. Differentiation using
$P_d'Q_d-P_dQ_d'=A_d^2$ gives
$D_{\R^2}\Psi_d=|\Psi_d|^2\Psi_d$ on $\R^2$. Away from the poles of
$R_d$, we have $|\Psi_d|^2=2|R_d'|/(1+|R_d|^2)$, so
$|\Psi_d|^4\,d\vol_{\ig_{\R^2}}$ is the pullback of the round area
element in stereographic coordinates. The Euclidean energy is therefore
\[
 \frac14\int_{\R^2}|\Psi_d|^4\,d\vol_{\ig_{\R^2}}
 =\frac14\int_{\mbs^2}R_d^*d\vol_{\ig_0}
 =\pi d,
\]
where $d=1,3,4,5,\dots$.

To write out the components, change the constant spinor frame by
$U=\diag(1,-i)$ and put $\widetilde\Psi_d=U\Psi_d$.
In this frame,
$D_{\R^2}(u,v)=(-2i\pa_zv,\,-2i\pa_{\bar z}u)$.
For $d=2k+1$, this recovers the angular-momentum solutions of
\cite[Theorem~7.1]{BF2020} with angular parameter $k$ and scale one, namely
\[
 \widetilde\Psi_{2k+1}(z)=
 \frac{\sqrt{2(2k+1)}}{1+|z|^{4k+2}}
 \begin{pmatrix}
 z^k\\
 i\,z^{k+1}|z|^{2k}
 \end{pmatrix},
 \qquad k\ge0.
\]
For $d=2k$, $k\ge2$, set $t=k/(k-1)$ and
\[
 F_k(z)=|z-1|^2+|z|^{4k-2}|z-t^2|^2.
\]
Then
\[
 \widetilde\Psi_{2k}(z)=
 \frac{\sqrt{2(2k-1)}}{F_k(z)}
 \begin{pmatrix}
 z^{k-1}(z-t)(\bar z-1)\\
 i\,z^k|z|^{2k-2}(z-t^2)(\bar z-t)
 \end{pmatrix}.
\]
Their energies are $(2k+1)\pi$ and $2k\pi$, respectively.
\end{remark}

For a generalized metric $\rho^2\ig_0$ with nonnegative, nonzero
$\rho\in L^2(\mbs^2)$, we use the positive eigenvalues of the
weighted equation $D\va=\la\rho\va$ in the sense of
\cite[Definition~1.1]{Martynyuk2026}. Recall that $\La_j$ was defined in
the introduction by taking the infimum of the normalized $j$-th
positive eigenvalue over smooth conformal metrics.
By \cite[Proposition~7.1]{Martynyuk2026}, the same value is obtained
by taking the infimum over generalized conformal metrics, so
\[
 \La_j=\inf_{\substack{\rho\in L^2(\mbs^2)\\
                        \rho\ge0,\ \rho\not\equiv0}}
       \la_j^+(\rho^2\ig_0)\,\norm{\rho}_2.
\]
Fix a nontrivial solution $\psi$ of \eqref{eq:two-dimensional-critical} with energy
$\pi d$ and set $\rho=|\psi|^2$. Then $1$ is a positive eigenvalue for
$\rho^2\ig_0$ and $\norm{\rho}_2^2=4\pi d$. For every integer $k\ge1$,
the definition of $\La_{2k}$ and \cite[Theorem~1.3]{Martynyuk2026}
give
\[
 \la_{2k}^+(\rho^2\ig_0)
 \ge\frac{\La_{2k}}{\norm{\rho}_2}
 =\sqrt{\frac{k}{d}}.
\]
By quaternionic symmetry, the positive eigenvalues occur in pairs.
If $2k$ is the last index at which the eigenvalue $1$ occurs, the
bound gives $1\ge\sqrt{k/d}$, hence $k\le d$. The metric $\rho^2\ig_0$
attains $\La_{2d}$ precisely when this last index is $2d$.

\begin{proposition}\label{prop:second-pair-nonattainment}
On $\mbs^2$, the common value $\La_3=\La_4=2\sqrt{2\pi}$ computed in
\cite[Theorem~1.3]{Martynyuk2026} is unattained among generalized
conformal metrics $\rho^2\ig_0$ with
$\rho\in L^2(\mbs^2)$, $\rho\ge0$, and $\rho\not\equiv0$.
\end{proposition}
\begin{proof}
Suppose that $\rho^2\ig_0$ attains $\La_3$ and normalize
$\norm{\rho}_2=1$. Set $\la=\la_3^+(\rho^2\ig_0)$.
The upper bound from two spheres in
\cite[Theorem~8.2]{Martynyuk2026} gives
$\la^2\le\La_1^2+\La_2^2=8\pi$. By
\cite[Lemma~7.2 and Remark~7.1]{Martynyuk2026}, $\rho$ is smooth,
and the corresponding weighted eigenspinors are smooth by elliptic
regularity. The quaternionic structure commutes with $D$ and
multiplication by $\rho$, so the complex multiplicity of $\la$
is even.

If this multiplicity is two, $\la$ occupies the third and fourth
positions. Applying \cite[Lemma~7.2]{Martynyuk2026} at the third
index gives $\rho=|\va|^2$ for a weighted eigenspinor $\va$.
Then $\psi=\sqrt{\la}\,\va$ solves \eqref{eq:standard} with
energy $\la^2/4\le2\pi$. Proposition~\ref{prop:two-dimensional-gap}
implies that $\psi$ is a ground state. Its induced metric
$|\psi|^4\ig_0=\la^2\rho^2\ig_0$ is round up to a conformal
diffeomorphism, with $1$ as its first and second positive eigenvalues.
Rescaling to $\rho^2\ig_0$ gives $\la$ in these two positions as well.
Thus the complex multiplicity of $\la$ is at least $4$, a contradiction.

If the multiplicity is at least $4$, choose a point $p\in\mbs^2$
with $\rho(p)>0$. Evaluation at $p$ is a complex-linear map from the
$\la$-eigenspace to the spinor fiber at $p$, which has complex
dimension $2$. By the rank-nullity theorem, its kernel has complex
dimension at least $4-2=2$. Choose a nonzero eigenspinor
$\va$ in this kernel, so that $\va(p)=0$. Write $f=\la\rho$, so
$f(p)>0$, and let $Z$ be the zero set of $\va$. The nodal
result in \cite{Baer1997} gives a finite set $Z$ of zeros of finite
order. As in Step~1 of the proof of
Proposition~\ref{prop:two-dimensional-gap}, the spinorial Weierstra\ss\
representation associates to $\va$
a smooth $\R^3$-valued one-form on $\mbs^2$. This form is closed
on $\mbs^2\setminus Z$, hence on all of $\mbs^2$. Since the sphere
is simply connected, the form integrates to a branched conformal
immersion. On $\mbs^2\setminus Z$, its induced metric is
$\widehat\ig=|\va|^4\ig_0$ and its signed mean curvature is
$H_\va=f/|\va|^2$.

Let $K_\va$ be its Gauss curvature on the regular set and let
$N$ be the sum of the vanishing orders of $\va$ at its zeros.
At a zero of order $s$, the leading term in a flat conformal
chart is annihilated by the flat Dirac operator and has norm
$c|z|^s$ for some $c>0$. The radial derivative of $\log|\va|$
is therefore $s/|z|+O(1)$. This gives a contribution of $4\pi s$
to the curvature integral, as in the proof of
\cite[Proposition~8.1]{Ammann2009}. Gauss-Bonnet and the
inequality $K_\va\le H_\va^2$ give
\[
 8\pi\ge\la^2
 =\int_{\mbs^2\setminus Z}H_\va^2\,d\vol_{\widehat\ig}
 \ge\int_{\mbs^2\setminus Z}K_\va\,d\vol_{\widehat\ig}
 =4\pi(1+N)\ge8\pi.
\]
Equality makes the immersion totally umbilic. Since $\mbs^2\setminus Z$ is
connected, so the Codazzi equation gives a constant $H_\va$ there.
By continuity, $f=H_\va|\va|^2$ also holds at $p$. This gives
$f(p)=0$, contrary to $\rho(p)>0$.

Thus $\La_3$ is unattained. By quaternionic symmetry,
$\la_3^+(\rho^2\ig_0)=\la_4^+(\rho^2\ig_0)$ for every generalized
conformal metric, so $\La_4$ is also unattained.
\end{proof}

\subsection{A lower bound in terms of analytic constants}\label{subsec:gap-bound}

For solutions of \eqref{eq:standard}, the deficit in the sharp
Sobolev inequality \eqref{eq:sharp-spinorial-sobolev} is small when
the energy is close to $E_m$. The stability estimate
\eqref{eq:stability} bounds the distance to the optimizer cone
$\co$. The equation also controls the amplitude of a nearby
optimizer. Thus, after dividing the solution by $a_0$ and applying
a conformal transformation, we can estimate the difference from a
unit Killing spinor in the norm defined below. Local uniqueness
forces the solution to be a ground state when this norm is smaller
than the local uniqueness radius. Quantifying this energy range gives a lower bound for the energy gap above $E_m$.

Recall that $\om_m=\vol_{\ig_0}(\mbs^m)$ is the volume of the unit
round sphere. We use the normalized measure
$d\mu=\om_m^{-1}d\vol_{\ig_0}$ and write
$\norm{\cdot}_{\hat q,d\mu}$ for the corresponding $L^{\hat q}$ norm.
Set
\begin{\equ}\label{eq:normalized-XY}
 \norm{v}_X^2=\frac2m\int_{\mbs^m}(|D|v,v)\,d\mu,
 \qquad
 \norm{v}_Y=\frac2m\norm{Dv}_{\hat q,d\mu}.
\end{\equ}
The first expression is understood as the quadratic form of $|D|$.
In particular, every $\xi\in\ck$ satisfies
$\norm{\xi}_X=\norm{\xi}_Y=1$. Choose $B_m\ge1$ such that
\begin{\equ}\label{eq:normalized-elliptic}
 \norm{v}_X\le B_m\norm{v}_Y.
\end{\equ}
Such a constant follows from Lemma~\ref{lem:elliptic}.

Let $\cg_1=a_0^{-1}\cg$. Proposition~\ref{prop:local} and compactness
of $\ck$ give a radius $\rho_m>0$ with the following property:
if $u$ is a weak solution of
\begin{\equ}\label{eq:normalized-equation}
 Du=\frac m2|u|^{q-2}u
\end{\equ}
and $\inf_{\xi\in\ck}\norm{u-\xi}_X<\rho_m$, then $u\in\cg_1$.
We fix any radius with this property, and retain the stability constant
$b_m$ from \eqref{eq:stability}.

\begin{proposition}\label{prop:safe-gap}
The conclusion of Theorem~\ref{thm:gap} holds with
\begin{\equ}\label{eq:conditional-kappa}
 \ka_m^{\mathrm{safe}}=
 \min\left\{\frac12,
       \frac{m b_m\rho_m^2}{8(B_m+1)^2},
       \frac{\hat q\rho_m}{2}\right\}.
\end{\equ}
In particular,
\[
 \Ga_m\ge\ka_m^{\mathrm{safe}},\qquad
 \ka_m\ge\ka_m^{\mathrm{safe}}.
\]
The estimates used here leave the numerical values of
$b_m,B_m,\rho_m$ undetermined.
\end{proposition}
\begin{proof}
Let $\psi=a_0u$ be a nontrivial solution of \eqref{eq:standard}.
Then $u$ solves \eqref{eq:normalized-equation}, and
\[
 \tau:=\frac{\ce(\psi)}{E_m}=\int_{\mbs^m}|u|^q\,d\mu\ge1,
 \qquad \norm{u}_Y=\tau^{1/\hat q}.
\]
Dividing \eqref{eq:stability}, applied to $u$, by
$(m/2)^2\om_m^{2/\hat q}$ yields
\begin{\equ}\label{eq:normalized-stability}
 \norm{u}_Y^2-\frac2m\int(Du,u)\,d\mu
 \ge b_m\inf_{v\in\co}\norm{u-v}_Y^2.
\end{\equ}
Thus the stability constant is unchanged, and the left-hand side on
the solution is
\[
 \de(\tau)=\tau^{2/\hat q}-\tau=\tau(\tau^{1/m}-1).
\]

Fix $\vr>0$ and choose an approximate optimizer
$v=\ct_F(c\xi)\in\co$, with $c>0$ and $\xi\in\ck$, such that
\[
 \norm{u-v}_Y\le\sqrt{\de(\tau)/b_m}+\vr.
\]
Apply the inverse conformal action and write $w=\ct_F^{-1}u$.
By conformal invariance of $Y$,
\[
 \norm{w-c\xi}_Y=\norm{u-v}_Y\le\sqrt{\de(\tau)/b_m}+\vr,
 \qquad \norm{w}_Y=\tau^{1/\hat q}.
\]
Since $\norm{\xi}_Y=1$ and $\tau\ge1$, the reverse triangle inequality
gives
\[
 |c-1|\le |c-\tau^{1/\hat q}|+\tau^{1/\hat q}-1
 \le\norm{w-c\xi}_Y+\tau^{1/\hat q}-1.
\]
Using \eqref{eq:normalized-elliptic} and $\norm{\xi}_X=1$, we obtain
\begin{align*}
 \norm{w-\xi}_X
 &\le B_m\norm{w-c\xi}_Y+|c-1|\\
 &\le(B_m+1)\left(\sqrt{\frac{\de(\tau)}{b_m}}+\vr\right)
       +\tau^{1/\hat q}-1.
\end{align*}
Since $\vr$ is arbitrary, this proves
\begin{\equ}\label{eq:quantitative-normalization}
 \inf_{F,\,\xi\in\ck}\norm{\ct_Fu-\xi}_X
 \le (B_m+1)\sqrt{\frac{\tau(\tau^{1/m}-1)}{b_m}}
       +\tau^{1/\hat q}-1.
\end{\equ}
Neither this inequality nor its use requires the optimizer infimum to
be attained.

Suppose $1<\tau<1+\ka_m^{\mathrm{safe}}$. In particular, $\tau<2$.
By concavity of the powers $1/m$ and $1/\hat q$,
\[
 \tau^{1/m}-1\le\frac{\tau-1}{m},\qquad
 \tau^{1/\hat q}-1\le\frac{\tau-1}{\hat q}.
\]
The right-hand side of \eqref{eq:quantitative-normalization} is
therefore strictly less than
\[
 (B_m+1)\sqrt{\frac{2\ka_m^{\mathrm{safe}}}{m b_m}}
       +\frac{\ka_m^{\mathrm{safe}}}{\hat q}
 \le\frac{\rho_m}{2}+\frac{\rho_m}{2}=\rho_m.
\]
There consequently exist an actual conformal action and a unit
Killing spinor for which the distance is less than $\rho_m$.
Local uniqueness implies $u\in\cg_1$, contradicting $\tau>1$.
This excludes the asserted interval above $E_m$.
\end{proof}

\begin{remark}\label{rem:constants-not-explicit}
Proposition~\ref{prop:normal-spectral-gap} supplies a uniform explicit
inverse bound for the normal Hessian. To obtain a numerical value of
$\rho_m$, one would also have to quantify the nonlinear remainder and
the local coordinates of the critical manifold. More importantly, the
global constant $b_m$ used here is not explicitly estimated in
\cite{WZ2026}. Formula~\eqref{eq:conditional-kappa} therefore does
not by itself give an explicit function of the dimension, a monotonicity
statement for the optimal energy gap, or a positive bound uniform over
all dimensions. Thus $(2m+3)^{-1}$ supplies the linear part of a
quantitative estimate of $\rho_m$, which enters the energy-gap bound
together with $b_m$ and $B_m$.
\end{remark}

\subsection{Explicit upper bounds from excited solutions}

For any excited solution $\psi$, \eqref{eq:least-excited-definition2}
gives
\[
 \Ga_m\le\frac{\ce(\psi)}{E_m}-1.
\]
We obtain an explicit upper bound in each dimension $m\ge3$ by
constructing an excited solution and computing its energy.

\begin{proposition}\label{prop:explicit-gap-upper}
The set of excited solutions is nonempty in every dimension $m\ge3$.
More precisely,
\begin{\equ}\label{eq:explicit-gap-upper}
 \Ga_3\le\frac{98}{27},\qquad
 \Ga_m\le
 \left(1+\frac{2}{m-1}\right)^{m-1}-1
 \quad\text{for }m\ge4.
\end{\equ}
\end{proposition}
\begin{proof}
We first give the three-dimensional example; the underlying
constant-length eigenspinor is also described in \cite{MS2014}.
Write $V^\flat=\ig_0(V,\cdot)$ for the metric dual of a vector field
$V$, and let $*$ denote the Hodge star on the oriented round
three-sphere.
Let $\sa\in\{1,-1\}$
be the scalar by which the real Clifford volume element acts on
$\Sigma\mbs^3$. Thus, for a vector field $V$,
\[
 (*V^\flat)\cdot=-\sa V\cdot.
\]
Choose a unit-length Killing spinor $\va$ satisfying
$\nabla_X\va=\tfrac12X\cdot\va$, and choose a unit Hopf field
$V$ with
\[
 dV^\flat=-2\sa *V^\flat,\qquad \diver V=0.
\]
Both handedness choices occur among the standard Hopf fields on the
round three-sphere. The Clifford Leibniz identity gives
\[
 D(V\cdot\va)
 =dV^\flat\cdot\va
   -(\diver V)\va
   +\frac12(3-2)V\cdot\va
 =\frac52 V\cdot\va.
\]
Consequently $\chi=V\cdot\va$ has constant length one, and
$\psi=\tfrac52\chi$ is a solution of $D\psi=|\psi|\psi$.
Its energy ratio is
\[
 \frac{\ce(\psi)}{E_3}=\left(\frac53\right)^3=\frac{125}{27}.
\]
This proves the first bound in \eqref{eq:explicit-gap-upper}.

For $m\ge4$, we shall construct a linear spinor-valued polynomial $P$
on $\R^m$ with
\begin{\equ}\label{eq:linear-monogenic}
 D_{\R^m}P=0,\qquad |P(x)|=|x|.
\end{\equ}
Here $(e_j)_{j=1}^m$ is the standard orthonormal basis of $\R^m$,
and all Clifford products below use $e_j\cdot e_j=-1$.
If $m$ is even, take
\[
 A=\diag(I_{m/2},-I_{m/2}),\qquad
 P(x)=A(x)\cdot s,
\]
where $s$ is any unit constant spinor. Since $A$ is symmetric and
trace-free,
\[
 D_{\R^m}P=\sum_j e_j\cdot A(e_j)\cdot s
             =-\tr(A)s=0.
\]
Orthogonality of $A$ gives $|P(x)|=|x|$.

If $m\ge5$ is odd, put $\cos\theta=-1/4$ and let $A$ act by the
same planar rotation $R_\theta$ through angle $\theta$ on
$\spann\{e_1,e_2\}$ and
$\spann\{e_3,e_4\}$. On the remaining $m-4$ coordinate
axes choose signs $+1$ and $-1$ whose sum is one. Thus $A$ is orthogonal
and $\tr A=4\cos\theta+1=0$. The commuting skew-Hermitian
endomorphisms $e_1\cdot e_2$ and $e_3\cdot e_4$ have a nonzero joint
eigenspace on which their eigenvalues are $i$ and $-i$, respectively.
Indeed, multiplication by $e_1$ and $e_3$ interchanges their two
respective eigenvalues, so all four joint eigenspaces have the same
positive dimension. Choose a unit $s$ in this joint eigenspace and again
set $P(x)=A(x)\cdot s$. Then
\[
 D_{\R^m}P
 =\bigl[-\tr A
       +2\sin\theta(e_1\cdot e_2+e_3\cdot e_4)\bigr]s=0,
\]
and $|P(x)|=|x|$. This proves \eqref{eq:linear-monogenic} in every
required dimension, without enlarging the spinor representation.

Let $\varrho=|x|$, and set
\[
 \bt=1+\frac2{m-1},\qquad B=(m\bt)^{(m-1)/2}.
\]
Consider the spinor
\begin{\equ}\label{eq:monogenic-excited}
 \Psi(x)=B(1+\varrho^{2\bt})^{-m/2}
       \bigl[P(x)-\varrho^{\bt-1}x\cdot P(x)\bigr].
\end{\equ}
To check the equation directly, let $f$ and $g$ be real radial
coefficient functions, with primes denoting differentiation in
$\varrho$. The homogeneity and monogenicity of $P$ give
\begin{align*}
 D_{\R^m}(f(\varrho)P+g(\varrho)x\cdot P)
 &=-\bigl[\varrho g'+(m+2)g\bigr]P
       +\frac{f'}\varrho x\cdot P,\\
 |fP+gx\cdot P|^2&=\varrho^2(f^2+\varrho^2g^2).
\end{align*}
Here the mixed real scalar product vanishes because Clifford
multiplication by $x$ is skew-Hermitian. Substitute
\[
 f=B(1+\varrho^{2\bt})^{-m/2},\qquad
 g=-B\varrho^{\bt-1}(1+\varrho^{2\bt})^{-m/2}.
\]
Using $B^{2/(m-1)}=m\bt$ and
$m\bt=m+2+\bt-1$, the two coefficients in the preceding Dirac
identity give
\[
 D_{\R^m}\Psi=|\Psi|^{2/(m-1)}\Psi.
\]
The spinor extends across $x=0$ with bounded first derivatives, is
smooth away from zero, and satisfies
$|\Psi(x)|=O(|x|^{-m})$ as $|x|\to\infty$. In particular,
$\Psi\in H^1(\R^m)\subset\dot H^{1/2}(\R^m)$ and the equation holds
weakly everywhere. Under stereographic compactification the spinor
also extends as a weak solution at the added point. One can see this
directly in the inverted chart: up to the fixed spin isometry needed
to orient inversion, its expression is
\[
 B(1+|y|^{2\bt})^{-m/2}
       \bigl[P(y)+|y|^{\bt-1}y\cdot P(y)\bigr].
\]
It has bounded first derivatives and vanishes at $y=0$; in particular,
the boundary term in the weak equation tends to zero there. We do not
assert smoothness at these zeros, since the nonlinearity has a
nonintegral exponent in general.

Finally, \eqref{eq:monogenic-excited} has norm
\[
 |\Psi(x)|=B\varrho(1+\varrho^{2\bt})^{-(m-1)/2}.
\]
Writing $\om_{m-1}=\vol(\mbs^{m-1})$ and substituting
$t=\varrho^\bt$, we obtain
\[
 \int_{\R^m}|\Psi|^q\,d\vol_{\ig_{\R^m}}
 =\frac{(m\bt)^m}{\bt}\om_{m-1}
       \int_0^\infty\frac{t^{m-1}}{(1+t^2)^m}\,dt=\bt^{m-1}C_m.
\]
The last equality is the stereographic formula for $\om_m$.
Denote the stereographic compactification of $\Psi$ by $\widehat\Psi$.
This nontrivial solution on $\mbs^m$ satisfies
$\ce(\widehat\Psi)/E_m=\bt^{m-1}>1$, proving the remaining bounds.
\end{proof}

\begin{remark}\label{rem:upper-not-monotone}
The explicit upper-bound function
$x\mapsto(1+2/x)^x-1$ increases to $e^2-1$, since
\[
 \frac{d}{dx}\left[x\log\left(1+\frac2x\right)\right]
 =\log\left(1+\frac2x\right)-\frac2{x+2}>0.
\]
This is a statement about the examples' upper bounds. It implies neither
monotonicity of $\Ga_m$ nor a positive lower bound uniform in $m$.
\end{remark}

\section{Existence for prescribed coefficients}\label{sec:existence}

\subsection{The all-dimensional variational input}

We record carefully which part of the argument is taken from our previous work
\cite{Xu2024}. The final existence theorem there is two-dimensional,
but its spectral reduction, test spinors, barycenter construction and
subcritical min-max construction are carried out for $m\ge2$.
The only step that we replace is the final exclusion of concentration.

Fix a positive coefficient $H\in C^2(\mbs^m)$ and a level
$d\in(0,H_{\max})$. For $2<p\le q$, define
\[
 \cl_p(\psi)=\frac12\bigl(\norm{\psi^+}_{\msh}^2
                    -\norm{\psi^-}_{\msh}^2\bigr)
       -\frac1p\int_{\mbs^m}H|\psi|^p\,d\vol_{\ig_0},
\]
and put
\begin{\equ}\label{eq:tau-levels}
 \tau_H=\frac{E_m}{H_{\max}^{m-1}},\qquad
 \tau_d=\frac{E_m}{d^{m-1}}.
\end{\equ}
For each $u\in\msh^+$, let $h_p(u)\in\msh^-$ be the
unique maximizer of $v\mapsto\cl_p(u+v)$ on $\msh^-$,
as in \cite[Proposition~2.2]{Xu2024}. This defines the negative
spectral reduction $h_p:\msh^+\to\msh^-$.
The reduced functional is $I_p(u)=\cl_p(u+h_p(u))$. Let
\[
 \cn_p=\{u\in\msh^+\setminus\{0\}:I_p'(u)[u]=0\}
\]
be the corresponding Nehari-Pankov constraint.

When the coefficient varies at the critical exponent, we use a
coefficient subscript: for a bounded real function $K$, set
\[
 \cl_K(v)=\frac12\int_{\mbs^m}(Dv,v)\,d\vol_{\ig_0}
                -\frac1q\int_{\mbs^m}K|v|^q\,d\vol_{\ig_0}.
\]
Thus $\cl_H=\cl_q$ for the fixed coefficient $H$ above;
the subscript in $\cl_p$ instead denotes the exponent.

\begin{lemma}[An estimate independent of the minimum of the coefficient]
\label{lem:weighted-unweighted-bound}
Let $0\le H\le M$ almost everywhere on $\mbs^m$, where $M>0$.
Suppose $\psi\in\msh$ is a weak solution of
$D\psi=H|\psi|^{q-2}\psi$, and put
$A_H(\psi)=\int_{\mbs^m} H|\psi|^q\,d\vol_{\ig_0}$. There is a constant $K(m)>0$,
independent of $H$ and its minimum, such that
\begin{\equ}\label{eq:weighted-unweighted-bound}
 \norm{\psi}_{\msh}+\norm{\psi}_q
 \le K(m)M^{(\hat q-1)/\hat q}A_H(\psi)^{1/\hat q},
 \qquad \hat q=\frac{2m}{m+1}.
\end{\equ}
Recall from \eqref{eq:constants} that
$C_m=(m/2)^m\vol_{\ig_0}(\mbs^m)$, if $\psi\ne0$, then
\begin{\equ}\label{eq:nonnegative-weighted-lower}
 A_H(\psi)\ge C_mM^{-(m-1)}.
\end{\equ}
\end{lemma}
\begin{proof}
Since $\hat q(q-1)=q$ and $H^{\hat q}\le M^{\hat q-1}H$,
\[
 \norm{D\psi}_{\hat q}^{\hat q}=\int H^{\hat q}|\psi|^q\,d\vol_{\ig_0}
 \le M^{\hat q-1}A_H(\psi).
\]
Invertibility of $D$ and \eqref{eq:elliptic}, followed by the critical
Sobolev embedding, prove \eqref{eq:weighted-unweighted-bound}.
For a nonzero solution, $A_H(\psi)>0$, since otherwise $D\psi=0$.
The sharp inequality \eqref{eq:sharp-spinorial-sobolev} gives
\[
 \frac m2\om_m^{1/m}A_H(\psi)
 \le\norm{D\psi}_{\hat q}^2
 \le M^{2(\hat q-1)/\hat q}A_H(\psi)^{2/\hat q}.
\]
Using $2/\hat q-1=1/m$ and $2(\hat q-1)/\hat q=(m-1)/m$ gives
\eqref{eq:nonnegative-weighted-lower}.
\end{proof}

\begin{remark}\label{rem:nonnegative-limit-ps}
In particular, if $H_j\to H_0$ uniformly, $0\le H_j\le M$, and
$\psi_j$ are nontrivial critical solutions with bounded weighted critical
mass, then $(\psi_j)$ is bounded in $\msh$, even when
$\min H_j\to0$. It is also a Palais-Smale sequence for the limiting
functional $\cl_{H_0}$, after selecting a convergent energy
subsequence. Indeed,
\[
 \norm{(H_j-H_0)|\psi_j|^{q-2}\psi_j}_{\hat q}
 \le\norm{H_j-H_0}_{\infty}\norm{\psi_j}_q^{q-1}\to0,
\]
and
$|\cl_{H_0}(\psi_j)-\cl_{H_j}(\psi_j)|\to0$.
This observation supplies the required boundedness and approximate
criticality for studying a coefficient limit with zeros; nontriviality
and compactness of that limit require a separate argument.
\end{remark}

\begin{proposition}[Variational input from \cite{Xu2024}]\label{prop:minmax}
Let $m\ge2$ and let $H$ satisfy \textup{($\mathrm{H}^*$)}.
There are constants $\theta>0$ and $\widehat a<\tau_d$, a strictly
subcritical sequence $2<p_j<q$, $p_j\uparrow q$, and critical points
$\psi_j\in\msh$ of $\cl_{p_j}$ whose values
$c_j:=\cl_{p_j}(\psi_j)$ satisfy
\[
 \widehat a<2\tau_H-\theta,
 \qquad
 \tau_H+\theta\le c_j\le\widehat a.
\]
The level $c_j$ is the min-max value over continuous homotopies of a
map $\va_j:\ch\to\cn_{p_j}$ to a constant map, with
\[
 \max_{x\in\ch} I_{p_j}(\va_j(x))\le\tau_H+\theta/2.
\]
In particular, the critical points satisfy
\begin{\equ}\label{eq:subcritical}
 D\psi_j=H|\psi_j|^{p_j-2}\psi_j.
\end{\equ}
\end{proposition}
\begin{proof}
This is the conclusion of Steps~1-3 in the proof of
\cite[Theorem~1.3, pp.~234-237]{Xu2024}, before the two-dimensional
argument in Step~4. We recall the construction to specify its scope.

\smallskip\noindent\textit{Step 1: A test map below $\tau_d$.}
Choose a small $\sa>0$ and a small neighborhood $U$ of
$\Crit(H)\cap H^{-1}(d)$, as in the cited construction, so that
\[
 C_0=\{H\ge d+\sa\}\cup U\subset\mbs^m
\]
contains $\ch$ and $\ch$ is contractible in $C_0$.
The bubble test spinors and the negative spectral reduction then give
a map $\Phi_j:C_0\to\cn_{p_j}$.
The estimates near the finitely many
critical points at level $d$ give a strict upper bound below $\tau_d$.
Writing $\vr>0$ for the bubble scale, the error estimates are
stated separately for $m=2$ and $m\ge3$ in
\cite[Lemmas~4.1-4.2]{Xu2024}; the strict improvement is of order
$\vr^2|\log\vr|$ in dimension two and $\vr^2$ for $m\ge3$.
Since $d>2^{-1/(m-1)}H_{\max}$, the upper bound can be chosen below
$2\tau_H$ by a fixed amount.

\smallskip\noindent\textit{Step 2: Localization of the lowest sublevels.}
Fix $\theta\in(0,2\tau_H-\tau_d)$ small enough that the barycenter
map sends the reduced
sublevel
\[
 \{u\in\cn_{p_j}:I_{p_j}(u)\le\tau_H+\theta\}
\]
into $\ch_\de$ for all sufficiently large $j$, where $\de$ is
chosen as in \textup{($\mathrm{H}^*$)}. Its composition with $\Phi_j$ on $\ch$ is
homotopic to the inclusion $\ch\hookrightarrow\ch_\de$.
This is the low-level localization in Step~2 of the cited proof. It
uses the concentration alternative and the fact that a near-minimizing
sequence concentrates at the maximum set. It does not use the
geometric argument on surfaces. The localization holds for one fixed
positive level increment and all sufficiently large $j$. Taking the
bubble scale sufficiently small, and choosing the constants as in the
cited construction, we may also arrange that
\[
 \sup_{x\in\ch}I_{p_j}(\Phi_j(x))<\tau_H+\theta/2.
\]
Thus the prescribed initial family lies strictly below the level used
in the deformation argument.

\smallskip\noindent\textit{Step 3: The min-max level.}
Set $\va_j=\Phi_j|_{\ch}$ and consider the continuous homotopies
in $\cn_{p_j}$ from $\va_j$ to a constant map. This class is nonempty
by the contraction in $C_0$. If one such homotopy lay in the lowest
sublevel, its barycenter would contract $\ch$ in $\ch_\de$,
contrary to \textup{($\mathrm{H}^*$)}. The subcritical
Palais-Smale condition and the deformation argument therefore give a
critical level between $\tau_H+\theta$ and a fixed
$\widehat a<\tau_d$. Our choice of $\theta$ gives
$\widehat a<2\tau_H-\theta$. The lift $u\mapsto u+h_{p_j}(u)$ gives
the critical point of
$\cl_{p_j}$. None of these three steps uses $m=2$.
\end{proof}

The compactness input is the all-dimensional concentration analysis of
\cite[Propositions~2.7-2.8 and their proofs]{Xu2024}. Although the
statements there are formulated for a strictly increasing subcritical
sequence, the argument includes the critical endpoint, as recorded
below.

\begin{proposition}\label{prop:alternative}
Let $H\in C^2(\mbs^m)$ be positive and let
$p_j\in(2,q]$, $p_j\to q$. Suppose that $\psi_j$ are nontrivial solutions of
\eqref{eq:subcritical} with
\[
 \tau_H\le\cl_{p_j}(\psi_j)\le2\tau_H-\theta
\]
for a fixed $\theta>0$. After passing to a subsequence, either
$\psi_j$ converges strongly in $\msh$ to a nontrivial solution
of \eqref{eq:prescribed}, or $\psi_j\weakto0$ and has a single-bubble
decomposition with a strong $\msh$ remainder. In the latter
case there are a point
$a\in\Crit(H)$ and a nontrivial finite-energy spinor $\va$ on
$\R^m$ satisfying
\begin{\equ}\label{eq:limit-bubble}
 D_{\R^m}\va=H(a)|\va|^{q-2}\va
\end{\equ}
such that, writing $c=\lim_j\cl_{p_j}(\psi_j)$,
\begin{\equ}\label{eq:bubble-energy}
 c=\frac{H(a)}{2m}\int_{\R^m}|\va|^q\,d\vol_{\ig_{\R^m}},
 \qquad
 c\ge\frac{C_m}{2mH(a)^{m-1}}.
\end{\equ}
\end{proposition}
\begin{proof}
If infinitely many $p_j$ equal $q$, pass to that subsequence.
The extraction and energy-splitting arguments of
\cite[Propositions~2.7-2.8]{Xu2024} give strong convergence or a
single bubble with a strong $\msh$ remainder. In the bubbling
case, let $R_j\to0$ denote the concentration scales used in the
rescaling. Since $p_j=q$, the scaling factor is identically
\[
 R_j^{(m-1)(q-p_j)/2}=1.
\]
Write $\de_a$ for the unit Dirac measure at $a$.
Then $|\psi_j|^q\,d\vol_{\ig_0}$ converges to a positive multiple of
$\de_a$. Applying the Kazdan-Warner identity
\eqref{eq:kw-obstruction} gives $X(H)(a)=0$ for every conformal
vector field $X$. Such fields span the tangent space at $a$, hence
$a\in\Crit(H)$.

Otherwise, eventually $p_j<q$, and a strictly increasing subsequence
converges to $q$. The conclusions for nontrivial solution sequences in
\cite[Proposition~2.8 and Lemmas~5.3-5.4]{Xu2024} then give
$a\in\Crit(H)$ and rescaling defect $1$. In both cases the strong
remainder and the identity
\[
 \cl_{p_j}(\psi_j)
 =\left(\frac12-\frac1{p_j}\right)
   \int H|\psi_j|^{p_j}\,d\vol_{\ig_0}
\]
give the energy identity as in \cite[(4.24)]{Xu2024}.
The lower bound follows by applying the classical ground-state bound
to $H(a)^{(m-1)/2}\va$.
\end{proof}

\begin{remark}\label{rem:critical-endpoint-defect}
For general Palais-Smale sequences with $p_j<q$ and concentration
scales $R_j\to0$, the bubble equation may contain the defect
$\la=\lim_jR_j^{(m-1)(q-p_j)/2}\le1$.
For sequences of nontrivial solutions, \cite{Xu2024} proves
$\la=1$ and $a\in\Crit(H)$. If $p_j=q$, the identity $\la=1$ is automatic,
even for approximate critical sequences; localization in $\Crit(H)$
still requires an additional condition, such as requiring each $\psi_j$ to solve its equation.
\end{remark}

\subsection{Compactness in the energy window}\label{subsec:energy-window}

\begin{proposition}\label{prop:exact-endpoint-compactness}
Let $H\in C^2(\mbs^m)$ be positive, and suppose a closed interval
$[\ell,u]$ satisfies
\begin{\equ}\label{eq:exact-compactness-window}
 \tau_H<\ell\le u<(1+\ka_m)\tau_H,\qquad
 \frac{E_m}{H(a)^{m-1}}\notin[\ell,u]\quad\text{for every }a\in\Crit(H).
\end{\equ}
Let $p_j\in(2,q]$, $p_j\to q$, and let $\psi_j$ be nontrivial
solutions of $D\psi_j=H|\psi_j|^{p_j-2}\psi_j$ with
$\cl_{p_j}(\psi_j)\in[\ell,u]$. Then $\{\psi_j\}$ is relatively
compact in $\msh$. Every limit along a subsequence is a nontrivial
solution of \eqref{eq:prescribed} with weighted energy in $[\ell,u]$.
\end{proposition}
\begin{proof}
Pass to a subsequence with energies converging to $c\in[\ell,u]$.
Since $\ka_m\le1$, Proposition~\ref{prop:alternative} gives strong
convergence or a single bubble $\va$ at $a\in\Crit(H)$. In the latter
case $\Phi=H(a)^{(m-1)/2}\va$ solves the constant-coefficient equation
on $\R^m$. By \eqref{eq:bubble-energy},
\[
 \ce_{\R^m}(\Phi)=H(a)^{m-1}c
 \le H_{\max}^{m-1}u<(1+\ka_m)E_m.
\]
Corollary~\ref{cor:euclidean-gap} therefore gives
$\ce_{\R^m}(\Phi)=E_m$, and hence $c=E_m/H(a)^{m-1}$. This contradicts
\eqref{eq:exact-compactness-window}. The strongly convergent
subsequence has a critical-exponent limit of energy $c$ by
Proposition~\ref{prop:alternative}. The argument applies to every
subsequence.
\end{proof}

\begin{remark}
If $H$ has no critical values in $(d,H_{\max})$, with
$0<d<H_{\max}$, then any closed subinterval of
$(\tau_H,\min\{\tau_d,(1+\ka_m)\tau_H\})$ satisfies
\eqref{eq:exact-compactness-window}. Indeed, the ground-state bubble
energy at a critical point is either $\tau_H$ or at least $\tau_d$.
The localization to critical points uses the equations satisfied by
the $\psi_j$. The same conclusion need not hold for general
Palais-Smale sequences.
\end{remark}

\begin{corollary}[$C^1$ perturbations inside the energy window]
\label{cor:C1-window-compactness}
Let $H\in C^2(\mbs^m)$ be positive, and let $[\ell,u]$ satisfy
\eqref{eq:exact-compactness-window}.
Suppose that positive $H_j\in C^2(\mbs^m)$ converge to $H$ in $C^1$,
and that $\psi_j$ are nontrivial solutions of
$D\psi_j=H_j|\psi_j|^{q-2}\psi_j$ with
\[
 \ell\le \frac1{2m}\int_{\mbs^m}H_j|\psi_j|^q\,d\vol_{\ig_0}\le u.
\]
Then $\{\psi_j\}$ is relatively compact in $\msh$ and
uniformly bounded in $L^\infty$. Every limit is a nontrivial
solution with coefficient $H$.
\end{corollary}
\begin{proof}
Write the displayed energies as $c_j$, and take a subsequence with
$c_j\to c\in[\ell,u]$. Uniform boundedness of $H_j$ and the equations
satisfied by $\psi_j$ give
\[
 \norm{D\psi_j}_{\hat q}^{\,\hat q}
 \le(\max H_j)^{\hat q-1}2mc_j.
\]
Thus $\psi_j$ is bounded in $\msh$. Since $H_j\to H$
uniformly, it is a Palais-Smale sequence for the fixed critical
functional $\cl_q$, at level $c$.
The critical Palais-Smale concentration alternative, as in the
proofs of \cite[Propositions~2.7-2.8]{Xu2024}, gives strong
convergence or a single bubble at some point $a$.
In the latter case, $|\psi_j|^q\,d\vol_{\ig_0}$ converges to
$A\de_a$, with $A>0$.

Although a general Palais-Smale sequence need not concentrate at a
critical point, the equations with coefficients $H_j$ give the
Kazdan-Warner identity
\[
 0=\int X(H_j)|\psi_j|^q\,d\vol_{\ig_0}\to A\,X(H)(a)
\]
for every conformal vector field $X$. Here the limit uses
$C^1$ convergence of the coefficients. Such fields span the tangent
space at $a$, so $a\in\Crit(H)$.
The same energy argument as in
Proposition~\ref{prop:exact-endpoint-compactness} now excludes the
bubble. Strong $L^q$ convergence and uniform boundedness of the
coefficients allow small-energy regularity on a finite cover of the
sphere, giving a uniform $L^\infty$ bound. The reasoning applies
to every subsequence and proves relative compactness.
\end{proof}

\begin{remark}
The fixed lower separation $\ell>\tau_H$ is essential for this
argument. The concentration sequence in
Proposition~\ref{prop:weaker-convergence-blowup} approaches the
lower endpoint $\tau_H$, so it does not contradict
Corollary~\ref{cor:C1-window-compactness}.
\end{remark}

\begin{proof}[Proof of Theorem~\ref{thm:existence}]
Proposition~\ref{prop:minmax} gives nontrivial solutions with
$2<p_j<q$, $p_j\uparrow q$, and
\[
 \tau_H+\theta\le\cl_{p_j}(\psi_j)\le\widehat a<\tau_d.
\]
By \eqref{eq:pinching},
$\tau_d<(1+\ka_m)\tau_H$. The critical-value condition in
\textup{($\mathrm{H}^*$)} then gives \eqref{eq:exact-compactness-window}
with $\ell=\tau_H+\theta$ and $u=\widehat a$.
A subsequence converges strongly in $\msh$ to a nontrivial
critical solution $\psi$, whose energy belongs to
$[\tau_H+\theta,\widehat a]$.
Testing its equation gives
$\cl_q(\psi)=(2m)^{-1}\int_{\mbs^m} H|\psi|^q\,d\vol_{\ig_0}$.
The strict bounds $\tau_H<\cl_q(\psi)<\tau_d$ prove
\eqref{eq:existence-energy}.
\end{proof}

\begin{remark}[Why the critical Palais-Smale condition is not automatic]
\label{rem:critical-ps-window}
The endpoint inclusion does not by itself permit a direct critical
min-max argument. In fact, under \textup{($\mathrm{H}^*$)}, the functional
$\cl_q$ has a noncompact Palais-Smale sequence at every level
$c\in(\tau_H,\tau_d)$.
To see this, set $h=(E_m/c)^{1/(m-1)}\in(d,H_{\max})$ and choose
$a$ with $H(a)=h$. This point is noncritical by the
assumption. Let $B_\vr$ be conformal ground-state bubbles for the
constant coefficient $h$, concentrating at $a$. They satisfy
\[
 DB_\vr=h|B_\vr|^{q-2}B_\vr,\qquad
 |B_\vr|^q\,d\vol_{\ig_0}\weakto C_mh^{-m}\de_a,
 \qquad B_\vr\weakto0\text{ in }\msh.
\]
Their fixed positive critical mass excludes strong convergence, while
\[
 \|(h-H)|B_\vr|^{q-2}B_\vr\|_{\hat q}^{\hat q}
 =\int|H-h|^{\hat q}|B_\vr|^q\,d\vol_{\ig_0}\to0.
\]
Thus $\cl_q'(B_\vr)\to0$ in $\msh^*$ and
$\cl_q(B_\vr)\to E_m/h^{m-1}=c$.
These sequences need not concentrate in $\Crit(H)$.
This separates the compactness of nontrivial critical solutions from the
Palais-Smale condition used in a direct variational construction.
The subcritical problem supplies nontrivial solutions before we apply
Proposition~\ref{prop:exact-endpoint-compactness}; a direct critical
construction would need an additional argument controlling escape
through such approximate bubbles.
\end{remark}

\subsection{An explicit class of admissible coefficients}\label{subsec:examples}

We next show that the combined topological and pinching conditions
are nonempty in every dimension. The example also allows maximum sets
of positive dimension.

\begin{proposition}\label{prop:example}
Fix $m\ge3$ and $a\in(0,1)$. For all sufficiently small $\vr>0$,
the function
\begin{\equ}\label{eq:example}
 H(x)=1-\vr(x_{m+1}-a)^2,\qquad x\in\mbs^m,
\end{\equ}
satisfies \textup{($\mathrm{H}^*$)} and \eqref{eq:pinching}
with
\[
 d=1-\vr(1-a)^2.
\]
Consequently, \eqref{eq:prescribed} has a nontrivial solution for this
coefficient.
\end{proposition}
\begin{proof}
Put $t=x_{m+1}$, and denote the north and south poles by $N$ and $S$.
The maximum value is $H_{\max}=1$, attained on
\[
 \ch=\{t=a\}\cong\mbs^{m-1}.
\]
For small $\de$, the tubular neighborhood $\ch_\de$ retracts
onto $\ch$. Since $\mbs^{m-1}$ is not contractible, its inclusion
into that neighborhood is not null-homotopic. On the other hand,
\[
 \{H\ge d\}=\{t\ge2a-1\}
\]
is a closed spherical cap, hence a topological $m$-ball. Thus $\ch$
is contractible in $\{H\ge d\}$.

The critical set consists of $\ch$ and the two poles. Moreover,
\[
 H(N)=d,\qquad H(S)=1-\vr(1+a)^2<d.
\]
There is no critical value in $(d,1)$, and the only critical point at
level $d$ is $N$. Since $\Hess t=-tg_0$,
\[
 \Hess H(N)=2\vr(1-a)\ig_0>0.
\]
It remains to choose $\vr$ so that positivity and pinching hold.
It is enough to require
\begin{\equ}\label{eq:example-epsilon}
 0<\vr<\min\left\{
 \frac1{(1+a)^2},
 \frac{1-(1+\ka_m)^{-1/(m-1)}}{(1-a)^2}
 \right\}.
\end{\equ}
Then $H>0$ and
$d>(1+\ka_m)^{-1/(m-1)}$. Since $\ka_m\le1$, this also implies
the original inequality $d>2^{-1/(m-1)}$ in
\textup{($\mathrm{H}^*$)}. Theorem~\ref{thm:existence} applies.
\end{proof}

\begin{remark}[The coefficient conditions in the two applications]
\label{rem:LT-coefficient-comparison}
For the preceding quadratic example,
\[
 \De_{\ig_0}H=-2\vr(1-a^2)<0\quad\text{on }\ch,
 \qquad\De_{\ig_0}=\diver\nabla.
\]
It therefore satisfies the coefficient condition of
\cite[Theorem~1.3]{LT2026}. That theorem gives uniform bounds for
critical solutions in its specified mass range, including small $C^2$
perturbations. Theorem~\ref{thm:existence} supplies nontrivial solutions
through the topological construction and the energy gap.

The existence hypotheses also allow flat maxima. For example, set
\[
 H(x)=1-\vr(x_{m+1}-a)^4,\qquad 0<a<1,
\]
and choose
\[
 0<\vr<\min\left\{(1+a)^{-4},\,
 \frac{1-(1+\ka_m)^{-1/(m-1)}}{(1-a)^4}\right\}.
\]
Its maximum set is $\{x_{m+1}=a\}$, and the level
$d=1-\vr(1-a)^4$ gives the same superlevel cap as in
Proposition~\ref{prop:example}. The only critical point at level $d$
is $N$, with $\Hess H(N)=4\vr(1-a)^3g_0>0$.
Thus \textup{($\mathrm{H}^*$)} and the pinching hold, while
$\Hess H=0$ on the maximum set. This coefficient does not meet the
strict Laplacian hypothesis of \cite[Theorem~1.3]{LT2026}, but
Theorem~\ref{thm:existence} still gives a nontrivial solution.

A priori compactness and nonemptiness have distinct roles here.
For the positive height coefficient $c+a x_{m+1}$, $c>a>0$, the
Laplacian is negative at its unique maximum, whereas the
Kazdan-Warner identity excludes every nontrivial solution.
The nonemptiness conclusion in our existence theorem therefore uses
the variational hypotheses in addition to compactness.
\end{remark}

\begin{corollary}[Coefficients with arbitrarily large contrast]
\label{cor:large-contrast}
Put $\theta_m=(1+\ka_m)^{-1/(m-1)}$ and fix $a\in(0,1)$ such that
\begin{\equ}\label{eq:large-contrast-a}
 \rho_a=\left(\frac{1-a}{1+a}\right)^2<1-\theta_m.
\end{\equ}
For every $0<\de<1$, the coefficient
\begin{\equ}\label{eq:large-contrast-coefficient}
 H_\de(x)=1-\frac{1-\de}{(1+a)^2}(x_{m+1}-a)^2
\end{\equ}
satisfies the hypotheses of Theorem~\ref{thm:existence}, with
\[
 (H_\de)_{\max}=1,\quad (H_\de)_{\min}=\de,\quad
 d_\de=1-(1-\de)\rho_a>1-\rho_a>\theta_m.
\]
The resulting solutions satisfy the uniform estimates
\begin{\equ}\label{eq:large-contrast-energy}
 C_m<\int H_\de|\psi_\de|^q\,d\vol_{\ig_0}
 <\frac{C_m}{(1-\rho_a)^{m-1}}<(1+\ka_m)C_m,
 \qquad \|\psi_\de\|_{\msh}\le C(m,a).
\end{\equ}
In particular, $(H_\de)_{\max}/(H_\de)_{\min}=1/\de$ is unbounded.
\end{corollary}
\begin{proof}
Take $\vr=(1-\de)/(1+a)^2$ in Proposition~\ref{prop:example}.
The sufficient restrictions in \eqref{eq:example-epsilon} follow from
$\de>0$ and $(1-\de)\rho_a<1-\theta_m$. The formulas for the
extrema and $d_\de$ are immediate, and
$d_\de-\de=(1-\de)(1-\rho_a)>0$.
Theorem~\ref{thm:existence} gives the mass bounds.
Lemma~\ref{lem:weighted-unweighted-bound} gives the uniform Hilbert bound.
\end{proof}

\begin{remark}\label{rem:large-contrast-nonperturbative}
The choice of $a$ is independent of $\de$. After normalizing
$\max H_\de=1$, we have
\[
 \inf_{c\in\R}\|H_\de-c\|_{\infty}
 =\frac{1-\de}{2}\to\frac12.
\]
Thus the criterion imposes no global closeness to a constant.
A comparison with any particular earlier existence theorem still
requires checking its hypotheses.
\end{remark}

\subsection{Density and its topology}\label{subsec:stronger-density}

\begin{proof}[Proof of Corollary~\ref{cor:stronger-density} for one solution]
We first give a local construction near the unique global maximum of
a positive smooth Morse function $H$. Write $M=\max H$. In a fixed
Morse coordinate neighborhood of its maximum, centered at $y=0$, we
have
\[
 H(y)=M-|y|^2.
\]
Choose $g\in C^\infty([0,\infty))$ such that
\[
 \begin{gathered}
  g(0)=-1,\qquad g(1)=0,\qquad g(u)=-u\quad(u\ge2),\\
  g'(u)>0\quad(0\le u<1),\qquad g'(u)<0\quad(u>1).
 \end{gathered}
\]
For example, take
$g(u)=-(1-\chi(u))u-\chi(u)(u-1)^2$, where $\chi$ is a smooth
nonincreasing cutoff equal to $1$ on $[0,5/4]$ and $0$ on $[2,\infty)$.
On the transition interval both profiles are decreasing and
$-(u-1)^2+u>0$, which verifies the required derivative sign.
For sufficiently small $\rho>0$,
replace $H$ in the coordinate ball $|y|<2\rho$ by
\begin{\equ}\label{eq:volcano-modification}
 H_\rho(y)=M+\rho^2g\bigl(|y|^2/\rho^2\bigr),
 \qquad d_\rho=M-\rho^2,
\end{\equ}
and set $H_\rho=H$ outside this ball. The matching is smooth because
$g(u)=-u$ for $u\ge2$. Positivity follows for small $\rho$ from the
uniform convergence proved below.

The maximum of $H_\rho$ is still $M$, and its maximum set is precisely
the coordinate sphere $|y|=\rho$. This sphere is not contractible in
its sufficiently small tubular neighborhood. The equation $g(u)=-1$
has precisely one solution $u_*>1$, apart from $u=0$, and $u_*<2$.
Since $H$ has a unique global maximum, shrinking $\rho$ also ensures
that all points outside the fixed Morse neighborhood lie below
$d_\rho$. Consequently
\[
 \{H_\rho\ge d_\rho\}=\{|y|\le\rho\sqrt{u_*}\},
\]
a closed ball containing the maximum sphere. There is no critical
value in $(d_\rho,M)$, and the only critical point at level $d_\rho$
is $y=0$. At this point the coordinate Hessian is
$2g'(0)\id>0$; as it is a critical point, this is also
the intrinsic Hessian quadratic form. Moreover $d_\rho>\min H_\rho$.
Thus the topological and critical-point requirements in
\textup{($\mathrm{H}^*$)} hold. Finally,
\[
 \left(\frac{M}{d_\rho}\right)^{m-1}\to1,
\]
so for small $\rho$ the additional pinching and the original
$2^{-1/(m-1)}M<d_\rho$ restriction both hold.

To estimate the modification, put $B(u)=g(u)+u$. In the fixed
coordinates,
\[
 f_\rho:=H_\rho-H
 =\rho^2 B\bigl(|y|^2/\rho^2\bigr),\qquad
 \supp f_\rho\subset\{|y|\le2\rho\}.
\]
Since $B$ is smooth and vanishes for $u\ge2$, differentiation gives
\begin{\equ}\label{eq:volcano-derivative-bounds}
 \|f_\rho\|_{C^0}\le C\rho^2,\qquad
 \|\nabla f_\rho\|_{C^0}\le C\rho,\qquad
 \|\nabla^2f_\rho\|_{C^0}\le C.
\end{\equ}
The constants may depend on $H$ and its fixed coordinate neighborhood,
but not on $\rho$. These bounds and the support volume estimate imply
\begin{\equ}\label{eq:volcano-strong-norms}
 \|H_\rho-H\|_{W^{2,s}}\le C_s\rho^{m/s}
 \quad(1\le s<\infty),\qquad
 \|H_\rho-H\|_{C^{1,\al}}
 \le C_\al\rho^{1-\al}
 \quad(0<\al<1).
\end{\equ}
For the second estimate, interpolate the uniform bound
$\|\nabla f_\rho\|_\infty=O(\rho)$ with its uniform Lipschitz
bound in the fixed charts. In particular,
$\|H_\rho-H\|_{C^1}=O(\rho)$.

It remains to approximate a general coefficient in any of the stated
spaces. Smooth functions are dense in $C^1$ and in $W^{2,s}$ for finite
$s$. In the last case, $s>m/2$ ensures uniform convergence as well,
so a sufficiently accurate
smooth approximation preserves a strictly positive lower bound.
An arbitrarily small smooth perturbation makes the approximation Morse
with distinct critical values, hence with a unique global maximum,
while preserving positivity. Its size can be made arbitrarily small
in each norm under consideration. Applying the construction above and taking $\rho$ sufficiently small
gives a coefficient in $\msd_m$ within any prescribed tolerance.
Theorem~\ref{thm:existence} gives a nontrivial solution with weighted
energy between $\tau_H$ and $\tau_d<(1+\ka_m)\tau_H$ for this coefficient.
\end{proof}

To compare min-max levels for different coefficients, write
$\cl_{p,K}$ for $\cl_p$ with a positive $C^2$ coefficient $K$.
On the unit sphere $S^+$ of $\msh^+$, set
\[
 J_{p,K}(w)=\max_{t\ge0,\,v\in\msh^-}\cl_{p,K}(tw+v).
\]
The radial identification of the Nehari-Pankov constraint with $S^+$
represents the constrained functional by $J_{p,K}$
\cite[Section~2.2, p.~217]{Xu2024}. In particular, a homotopy class
defined on $S^+$ can be kept fixed when the coefficient changes.

\begin{lemma}
\label{lem:subcritical-minmax-persistence}
Fix $p_0\in(2,q)$, let $p_j\in[p_0,q)$, and let
$\va_j:X\to S^+$ be continuous maps on a nonempty compact space $X$.
Suppose that each $\va_j$ is homotopic to a constant map. Write
$\cf_j$ for the class of these homotopies, and for a positive
coefficient $K\in C^2(\mbs^m)$ set
\[
 c_j(K)=\inf_{\gamma\in\cf_j}\max_{(x,t)\in X\times[0,1]}
              J_{p_j,K}(\gamma(x,t)).
\]
Suppose a positive $H\in C^2(\mbs^m)$ and fixed numbers
$b,\ell,u$ satisfy
\[
 \max_X J_{p_j,H}(\va_j)\le b<\ell\le c_j(H)\le u
 \qquad\text{for every }j.
\]
For $0<\eta<(\ell-b)/2$, there is $\vr\in(0,1)$, independent of $j$,
such that whenever $(1-\vr)H\le K\le(1+\vr)H$, we have
\[
 \max_X J_{p_j,K}(\va_j)\le b+\eta<\ell-\eta,\qquad
 c_j(K)\in[\ell-\eta,u+\eta].
\]
Each $c_j(K)$ is the energy of a nontrivial critical point of
$\cl_{p_j,K}$.
\end{lemma}
\begin{proof}
For a constant $\la>0$, scaling the spinor by
$\la^{-1/(p-2)}$ in the fiber defining $J_{p,K}$ gives
$J_{p,\la H}=\la^{-2/(p-2)}J_{p,H}$. Pointwise ordering in the
coefficient therefore gives
\[
 (1+\vr)^{-2/(p-2)}J_{p,H}
 \le J_{p,K}
 \le (1-\vr)^{-2/(p-2)}J_{p,H}.
\]
Use these inequalities with the fixed homotopy class $\cf_j$ and
with its initial map $\va_j$. Writing $\nu=2/(p_0-2)$, we obtain
\[
 (1+\vr)^{-\nu}\ell\le c_j(K)\le(1-\vr)^{-\nu}u,\qquad
 \max_X J_{p_j,K}(\va_j)\le(1-\vr)^{-\nu}b.
\]
Taking $\vr$ sufficiently small gives the asserted bounds for all $j$.
Since $p_j<q$, the constrained functional satisfies the Palais-Smale
condition. The deformation at level $c_j(K)$ fixes the initial map,
whose values lie strictly below the level, and sends constant maps
to constant maps. It follows as in Proposition~\ref{prop:minmax}
that $c_j(K)$ is a constrained critical value. The lift from the
Nehari-Pankov constraint gives the stated critical point.
\end{proof}

\begin{proof}[Proof of the multiplicity assertion in Corollary~\ref{cor:stronger-density}]
The conditions in \eqref{eq:exact-compactness-window} persist under
small enlargements of the interval and small $C^1$ changes of the
coefficient. The strict endpoint inequalities persist by continuity
of $\tau_H$. For the exclusion of bubble energies, observe that if
$K_j\to H$ in $C^1$ and $a_j\in\Crit(K_j)$ converge to $a$, then
$a\in\Crit(H)$ and $E_m/K_j(a_j)^{m-1}\to E_m/H(a)^{m-1}$.

Start with a smooth Morse approximation of the target coefficient
with a unique maximum. The modification \eqref{eq:volcano-modification}
and Proposition~\ref{prop:minmax} give the first family of
subcritical min-max values. Radial projection to $S^+$ gives the
initial maps for Lemma~\ref{lem:subcritical-minmax-persistence}.
Choose a closed interval inside $(\tau_H,\tau_d)$ containing the
min-max values and with its lower endpoint above the initial-map
energies. The bounds in Proposition~\ref{prop:minmax} allow both
separations to be uniform in $j$. This interval satisfies
\eqref{eq:exact-compactness-window} by \textup{($\mathrm{H}^*$)}
and \eqref{eq:pinching}.

At each subsequent step, first enlarge the retained intervals
slightly so that they remain disjoint and still satisfy
\eqref{eq:exact-compactness-window}. Choose the next perturbation
small enough in $C^0$ to preserve the min-max families in these
intervals by Lemma~\ref{lem:subcritical-minmax-persistence}, and
small enough in $C^1$ to preserve \eqref{eq:exact-compactness-window}.
To obtain the new family, approximate the current coefficient by a
smooth Morse function with distinct critical values and modify it
near its maximum by \eqref{eq:volcano-modification}. The Morse
approximation changes the maximum by an arbitrarily small amount.
As $\rho\to0$, the new level $d_\rho$ tends to that maximum.
Thus the upper energy bound $E_m/d_\rho^{m-1}$ for the new family
can be chosen below the lower endpoints of all the retained intervals.
The new interval again satisfies \eqref{eq:exact-compactness-window}
and has uniform separation from its initial-map energies.

After $k$ steps, fix the final coefficient $H\in\msd_m$ and
take the subcritical exponents to $q$ in each family.
Proposition~\ref{prop:exact-endpoint-compactness} yields one
nontrivial solution in each of the $k$ disjoint intervals.
Each modification can be made arbitrarily small in the norm used
for approximation by \eqref{eq:volcano-strong-norms}. Keeping the
sum of the changes below the prescribed tolerance gives the
asserted density.
\end{proof}

\begin{remark}
The same argument gives density in the positive cone of the 
H\"older space, defined as the closure of $C^\infty$ in $C^{1,\al}$.
It does not give density in the full space $C^{1,\al}$, where smooth
functions are not dense. Likewise, the estimate gives no convergence in
$W^{2,\infty}$ or in $C^2$.
\end{remark}

\begin{proof}[Proof of Theorem~\ref{thm:baire-multiplicity}]
\smallskip\noindent\textit{Step 1: Compactness for a $C^1$ limit coefficient.}
Let $K_j,K\in C^1_+(\mbs^m)$ with $K_j\to K$ uniformly, and write
$K_{\max}=\max_{\mbs^m}K$ and $\tau_K=E_m/K_{\max}^{m-1}$.
Suppose that $\psi_j$ are nontrivial solutions with coefficients $K_j$
and weighted energies $\cl_{K_j}(\psi_j)\to c\in(0,2\tau_K)$. By
Lemma~\ref{lem:weighted-unweighted-bound}, $\{\psi_j\}$ is bounded in
$\msh$. As in Remark~\ref{rem:nonnegative-limit-ps}, it is a
Palais-Smale sequence for $\cl_K$ at level $c$. After passing to a
subsequence, $\psi_j\weakto\psi$ in $\msh$, and the convergence is
strong in $L^2$ and holds almost everywhere. The nonlinearity is
bounded in $L^{\hat q}$ and converges almost everywhere, so $\psi$
solves the equation with coefficient $K$.

Write $v_j=\psi_j-\psi$. The Brezis-Lieb splitting for the critical
power and its derivative gives
\[
 Dv_j-K|v_j|^{q-2}v_j\to0\quad\text{in }\msh^*,\qquad
 \frac{1}{2m}\int_{\mbs^m}K|v_j|^q\,d\vol_{\ig_0}
   \to c-\cl_K(\psi).
\]
Set $w_j=K^{(m-1)/2}v_j$. Multiplication by a $C^1$ function is
bounded on $\msh$ by interpolation between $L^2$ and $H^1$, and on
$\msh^*$ by duality. The Dirac product rule gives
\[
 Dw_j-|w_j|^{q-2}w_j
 =K^{(m-1)/2}\bigl(Dv_j-K|v_j|^{q-2}v_j\bigr)
  +\nabla(K^{(m-1)/2})\cdot v_j\to0\quad\text{in }\msh^*.
\]
The last term converges to zero in $L^2$. Thus $w_j\weakto0$ in
$\msh$ and, along a subsequence, $\{w_j\}$ is a Palais-Smale sequence for
$\ce$ with
\[
 \lim_j\ce(w_j)=\lim_j\frac{1}{2m}\int_{\mbs^m}|w_j|^q\,d\vol_{\ig_0}
 \le K_{\max}^{m-1}\bigl(c-\cl_K(\psi)\bigr).
\]
If $\psi\ne0$, then $\cl_K(\psi)\ge\tau_K$ by
\eqref{eq:nonnegative-weighted-lower}. The energy level of $\{w_j\}$ is then
below $E_m$, so global compactness \cite[Theorem~5.2]{Isobe2011}
gives $w_j\to0$ in $\msh$. It follows that $\psi_j\to\psi$ in
$\msh$ and $\cl_K(\psi)=c$. If $\psi=0$, the energy level of $\{w_j\}$ is
positive and below $2E_m$. The same decomposition gives a single
bubble at some $a\in\mbs^m$ with a strong remainder. Its profile
$\Phi$ satisfies $D_{\R^m}\Phi=|\Phi|^{q-2}\Phi$. Since
$|\psi_j|^q=K^{-m}|w_j|^q$, we obtain
\begin{\equ}\label{eq:baire-bubble-energy}
 c=K(a)^{-(m-1)}\ce_{\R^m}(\Phi),\qquad
 |\psi_j|^q\,d\vol_{\ig_0}\weakto
    2m K(a)^{-m}\ce_{\R^m}(\Phi)\,\de_a.
\end{\equ}
The corresponding profile for $\psi_j$ is
$K(a)^{-(m-1)/2}\Phi$.

Now suppose that $K_j$ are smooth, $K_j\to K$ in $C^1$, and the
energies belong to a closed interval $[\ell,u]$ satisfying
\eqref{eq:exact-compactness-window} with $H$ replaced by $K$.
In the bubbling case, the Kazdan-Warner identity
\eqref{eq:kw-obstruction}, applied to $\psi_j$ and passed to the limit
using \eqref{eq:baire-bubble-energy}, gives $X(K)(a)=0$ for all
conformal vector fields $X$. Thus $\nabla K(a)=0$.
Since $K(a)^{m-1}c<(1+\ka_m)E_m$,
Corollary~\ref{cor:euclidean-gap} gives $\ce_{\R^m}(\Phi)=E_m$.
Then $c=E_m/K(a)^{m-1}\in[\ell,u]$, contrary to
\eqref{eq:exact-compactness-window}. Hence $\psi_j$ converges strongly
along a subsequence to a nontrivial solution with coefficient $K$ and
energy in $[\ell,u]$.

\smallskip\noindent\textit{Step 2: Open dense sets of coefficients.}
Fix $n\ge1$. Let $\mathcal U_n$ be the interior in $C^1_+(\mbs^m)$
of the set of coefficients $K$ for which \eqref{eq:prescribed} has a
nontrivial solution with weighted energy in
\[
 \bigl(\tau_K,\min\{\tau_K+1/n,(1+\ka_m)\tau_K\}\bigr).
\]
To see that $\mathcal U_n$ is dense, start with a smooth Morse
approximation with a unique global maximum to any positive $C^1$
coefficient. In the local modification \eqref{eq:volcano-modification}, the level $d_\rho$
approaches the maximum as $\rho\to0$. We may therefore take the
modified coefficient $H$ arbitrarily close to the target in $C^1$
and arrange that
\[
 \tau_{d_\rho}<\min\{\tau_H+1/n,(1+\ka_m)\tau_H\}.
\]
Proposition~\ref{prop:minmax} gives subcritical min-max values in a
closed interval strictly inside $(\tau_H,\tau_{d_\rho})$, with the
initial-map energies uniformly below that interval. By
Lemma~\ref{lem:subcritical-minmax-persistence}, these values remain in a
slightly enlarged closed interval $[\ell,u]$ for smooth $K$
sufficiently close to $H$. Shrinking the $C^1$ neighborhood of $H$
if necessary, the same interval satisfies
\eqref{eq:exact-compactness-window} for every $K$ in the neighborhood,
by the critical-value argument in the proof of
Corollary~\ref{cor:stronger-density}. We also retain $u<\tau_K+1/n$
by continuity of $\tau_K$. For smooth $K$,
Proposition~\ref{prop:exact-endpoint-compactness} then gives a
critical solution with energy in $[\ell,u]$. Approximation by smooth
coefficients in the same neighborhood and Step~1 give a solution in
that interval for every $C^1$ coefficient in the neighborhood.
Hence $H\in\mathcal U_n$.

\smallskip\noindent\textit{Step 3: The fixed-coefficient sequence.}
Since $C^1_+(\mbs^m)$ is open in the Banach space $C^1(\mbs^m)$, the
Baire category theorem shows that $\bigcap_{n\ge1}\mathcal U_n$ is
dense. For each $H$ in this intersection, the definition of
$\mathcal U_n$ gives solutions for the same $H$ with energies strictly
above $\tau_H$ and tending to $\tau_H$. Passing to a subsequence, we
may take these energies strictly decreasing. A constant coefficient
cannot belong to the intersection, by \eqref{eq:constant-scaling}
and the energy gap.

For a nonconstant coefficient $H$, the lower bound $\tau_H$ cannot
be attained by a nontrivial solution. Indeed, equality in the proof
of \eqref{eq:nonnegative-weighted-lower} implies $H=H_{\max}$ wherever
the spinor is nonzero, almost everywhere. The spinor then solves the
constant-coefficient equation with coefficient $H_{\max}$ and, after
rescaling, is a ground state. Since ground-state spinors are nowhere
zero, continuity of $H$ would give $H\equiv H_{\max}$. Apply Step~1
with $K_j=K=H$ and $c=\tau_H$. The strong convergence alternative is
excluded by this equality argument. Formula
\eqref{eq:baire-bubble-energy} and the ground-state lower bound give
\[
 \tau_H=H(a)^{-(m-1)}\ce_{\R^m}(\Phi)
       \ge \frac{E_m}{H(a)^{m-1}}\ge\tau_H.
\]
Thus $a$ is a global maximum of $H$ and $\Phi$ is a ground-state
bubble. Finally, any subsequence bounded in $L^\infty$ would, after
extraction, converge to zero in $L^2$ and hence in $L^q$, contrary to
its positive limiting weighted energy. This proves
$\norm{\psi_j}_\infty\to\infty$.
\end{proof}

\begin{proposition}[$C^2$-open nonexistence]
\label{prop:open-nonexistence}
Let $m\ge2$ and $H_0(x)=c+a x_{m+1}$, where $c>a>0$.
There exists $\eta>0$ such that every
$H\in C^2(\mbs^m)$ with
\[
 \|H-H_0\|_{C^2}<\eta
\]
is positive and admits no nontrivial weak solution of
\eqref{eq:prescribed}.
\end{proposition}
\begin{proof}
Put $t=x_{m+1}$, and let $N,S$ denote the north and south poles.
The critical points of $H_0$ are its nondegenerate maximum $N$ and
minimum $S$. Since $\Hess H_0=-atg_0$, a sufficiently small
$C^2$ perturbation has exactly two critical points $p_+,p_-$, near
$N,S$, respectively, with the same Hessian signs. Indeed, strict
concavity and convexity in fixed small geodesic balls give uniqueness
there, existence follows from the implicit function theorem, and the
gradient is bounded away from zero on their complement.

We normalize this pair by a conformal transformation close to the
identity. For completeness, identify both $T_N\mbs^m$ and
$T_S\mbs^m$ with $\R^m$. The conformal vector fields
$x\mapsto Bx$ with $B$ skew-symmetric, together with the gradients of
the first spherical harmonics, have surjective joint evaluation at
$N,S$. More explicitly, given $v,w\in\R^m$, take a rotation
field with values $(v-w)/2$ at $N$ and $-(v-w)/2$ at $S$, and the
gradient of the horizontal linear function with coefficient
$(v+w)/2$. Their sum has values $v,w$ at the two poles. The inverse
function theorem, applied to the time-one flows of this $2m$-parameter
family, therefore supplies a conformal transformation $A$ with
\[
 A(N)=p_+,\qquad A(S)=p_-,\qquad A\to\id
 \quad\text{as }H\to H_0\text{ in }C^2.
\]
In particular $K:=H\circ A$ is $C^2$-close to $H_0$ and is critical
at both poles.

We claim that $X_0(K)>0$ off $N,S$, where $X_0=\nabla t$.
Fix small geodesic balls around the poles. For $K$ sufficiently close
to $H_0$, there is $b>0$ such that
\[
 \Hess K\le-bg_0\text{ near }N,\qquad
 \Hess K\ge bg_0\text{ near }S.
\]
Along a unit-speed radial geodesic from $N$, the vanishing gradient at
$N$ gives $\pa_rK\le-br$. Since
$X_0=-\sin r\,\pa_r$ there, we obtain
$X_0(K)\ge br\sin r>0$ away from $N$. The same argument at $S$
uses $X_0=\sin r\,\pa_r$ and $\pa_rK\ge br$.
On the complement of these balls, positivity follows by $C^1$
closeness from
\[
 X_0(H_0)=a(1-t^2)>0.
\]
Consequently the conformal vector field $X=A_*X_0$ satisfies
\[
 X(H)>0\quad\text{on }\mbs^m\setminus\{p_+,p_-\}.
\]

The spinorial Kazdan-Warner identity
\eqref{eq:kw-obstruction}, proved in
\cite[Theorem~2 and Example~3(2)]{AHA2007}, applies to weak solutions
here as well. One can justify this directly by elliptic regularity:
a weak $H^{1/2}$ solution is bounded and then belongs to
$C^{2,\bt}$ for some $\bt>0$. Indeed, the map
$z\mapsto|z|^{q-2}z$ is $C^1$ with locally H\"older derivative,
and $H\in C^2$. Thus its conformal-flow orbit is differentiable in
$H^{1/2}$. Differentiating the critical action along this orbit gives
the identity, by conformal invariance of the quadratic term and the
critical exponent. It follows that
\[
 0=\int_{\mbs^m}X(H)|\psi|^q\,d\vol_{\ig_0}.
\]
The integrand is nonnegative and its coefficient is positive outside
two points. Hence $\psi=0$ almost everywhere. Finally, decrease
$\eta$ to preserve positivity of $H$.
\end{proof}

\begin{remark}
The density conclusions of Corollary~\ref{cor:stronger-density} are
compatible with Proposition~\ref{prop:open-nonexistence}: the local
modification changes the Hessian by a bounded amount on a shrinking
set, although its $W^{2,s}$ norm tends to zero for every finite $s$.
In particular solvability is not dense in the $C^2$ topology.
The smooth class $\msd_m$ is not dense in the $W^{2,\infty}$ or
$C^{1,1}$ norm either, even relative to positive smooth targets:
on smooth functions each of these norms controls the $C^2$ norm,
so the same nonexistence neighborhood excludes approximation to $H_0$.
\end{remark}

\subsection{Concentration under weaker convergence of the coefficients}
\label{subsec:weaker-coefficient-convergence}

We next construct nontrivial solutions exhibiting first-threshold
concentration as their smooth positive coefficients converge in
weaker topologies than $C^2$. The results of \cite{LT2026} describe
the local structure of these concentrating solutions. Here we obtain
such solutions by modifying a smooth coefficient near its unique
nondegenerate maximum.

\begin{proposition}\label{prop:weaker-convergence-blowup}
Let $m\ge3$ and
\[
 H_*(x)=c+a x_{m+1},\qquad c>a>0,\qquad M=c+a.
\]
There are positive smooth coefficients $H_j$ with maximum value $M$,
and nontrivial solutions
\[
 D\psi_j=H_j|\psi_j|^{q-2}\psi_j,
\]
such that
\begin{\equ}\label{eq:weaker-coefficient-limits}
 H_j\to H_*\quad\text{in }C^{1,\al}\text{ for every }0<\al<1,
 \quad
 H_j\to H_*\quad\text{in }W^{2,s}\text{ for every }1\le s<\infty.
\end{\equ}
Moreover, $\sup_j\norm{H_j}_{C^2}<\infty$, but $H_j$ does not
converge to $H_*$ in $C^2$. Let $p_N=(0,\ldots,0,1)$ denote the
north pole. Then
\begin{\equ}\label{eq:weak-coefficient-concentration}
 \psi_j\weakto0\text{ in }\msh,\qquad
 \norm{\psi_j}_{\infty}\to\infty,\qquad
 |\psi_j|^q\,d\vol_{\ig_0}\weakto\frac{C_m}{M^m}\de_{p_N}.
\end{\equ}
The convergence $\psi_j\to0$ is uniform on compact subsets of
$\mbs^m\setminus\{p_N\}$. After rescaling at the concentrating point
and passing to a subsequence, the solutions converge to a
constant-coefficient ground-state bubble. Their weighted masses satisfy
\[
 \int_{\mbs^m}H_j|\psi_j|^q\,d\vol_{\ig_0}\to\frac{C_m}{M^{m-1}}, \quad \text{as } j\to+\infty.
\]
\end{proposition}
\begin{proof}
\smallskip\noindent\textit{Step 1: The coefficients and nontrivial solutions.}
Take Morse coordinates $y$ near $p_N$ in which $H_*=M-|y|^2$.
Choose a smooth function $g:[0,\infty)\to\R$ such that
\[
 g(0)=-1,\quad g'(0)>0,\quad g(1)=0,\quad g''(1)<0,\quad
 g(s)=-s\quad(s\ge2),
\]
with $g'>0$ on $[0,1)$ and $g'<0$ on $(1,\infty)$.
For $\rho_j\searrow0$, replace $H_*$ in this coordinate
neighborhood by
\[
 H_j(y)=M+\rho_j^2 g(|y|^2/\rho_j^2),
\]
and retain $H_*$ outside it. The matching is smooth because $g=-s$
for $s\ge2$. As in the local modification above, the maximum set is
the sphere $|y|=\rho_j$, the center is a nondegenerate local minimum
at level $d_j=M-\rho_j^2$, and the corresponding superlevel is a
closed ball. There are no other critical values in $(d_j,M)$.
Thus \textup{($\mathrm{H}^*$)} holds. Since $M/d_j\to1$, the pinching
condition holds for all sufficiently large $j$.

The difference of the coefficients has the form
$\rho_j^2B(y/\rho_j)$ for a fixed compactly supported smooth
function $B$. The scaling estimates for this expression give
\eqref{eq:weaker-coefficient-limits} and the uniform $C^2$ bound.
The Hessian at $p_N$ changes from negative definite to positive
definite with a fixed nonzero size, so there is no $C^2$ convergence.
Discarding finitely many indices, all coefficients are positive and
bounded below by a fixed number.

Theorem~\ref{thm:existence} gives nontrivial solutions with
\begin{\equ}\label{eq:weak-convergence-weighted-squeeze}
 \frac{C_m}{M^{m-1}}
 <Q_j:=\int H_j|\psi_j|^q\,d\vol_{\ig_0}
 <\frac{C_m}{(M-\rho_j^2)^{m-1}}.
\end{\equ}
In particular, $Q_j\to Q_0:=C_m/M^{m-1}$.

\smallskip\noindent\textit{Step 2: A bounded Palais-Smale sequence.}
We use the coefficient-dependent critical functional $\cl_K$
defined at the beginning of this section.
Since $H_j\le M$ and $(q-1)\hat q=q$,
\[
 \norm{D\psi_j}_{\hat q}^{\,\hat q}
 =\int H_j^{\hat q}|\psi_j|^q\,d\vol_{\ig_0}
 \le M^{\hat q-1}Q_j.
\]
Lemma~\ref{lem:elliptic} gives a uniform $\msh$ bound.
Because $\norm{H_j-H_*}_{C^0}=O(\rho_j^2)$, the equations satisfied by $\psi_j$ yield
\begin{\equ}\label{eq:fixed-Hstar-PS}
 \cl_{H_*}(\psi_j)
 \to\frac{E_m}{M^{m-1}},\qquad
 \norm{\cl_{H_*}'(\psi_j)}_{\msh^*}
 =O(\rho_j^2).
\end{\equ}
Indeed, $\cl_{H_j}(\psi_j)=Q_j/(2m)$, and the differences
of the energies and derivatives are bounded by
$C\norm{H_j-H_*}_{C^0}$ using the uniform $L^q$ bound.

Every weak limit is a solution with coefficient $H_*$.
One may see this directly from compact convergence in $L^2$,
almost-everywhere convergence, and weak convergence of
$|\psi_j|^{q-2}\psi_j$ in $L^{\hat q}$.
The Kazdan-Warner identity \eqref{eq:kw-obstruction} excludes every
nontrivial such limit: for $X=\nabla x_{m+1}$,
\[
 X(H_*)=a(1-x_{m+1}^2)>0
\]
away from the poles. Thus $\psi_j\weakto0$ in $\msh$ and
$\psi_j\to0$ in $L^2$. A subsequence with bounded $L^\infty$ norm
would have $\int|\psi_j|^q\,d\vol_{\ig_0}\to0$, contradicting
\eqref{eq:weak-convergence-weighted-squeeze}. Hence
$\norm{\psi_j}_{\infty}\to\infty$.

\smallskip\noindent\textit{Step 3: The first bubble exhausts the mass.}
Standard blow-up rescaling for uniformly convergent positive
coefficients gives, along a subsequence, a point $p$ and a nontrivial
finite-mass Euclidean solution
\[
 D_{\R^m}\Phi=h|\Phi|^{q-2}\Phi,\qquad h=H_*(p).
\]
We use the rescaling construction of
\cite[Lemma~3.4 and Proposition~3.2]{LT2026}.
For this sequence the weighted energy squeeze identifies the limiting
bubble and its mass directly. Local convergence, conformal invariance
of the critical integral and the classical ground-state lower bound give
\begin{\equ}\label{eq:weighted-mass-exhaustion}
 Q_0\ge h\int_{\R^m}|\Phi|^q\,d\vol_{\ig_{\R^m}}
 \ge\frac{C_m}{h^{m-1}}
 \ge\frac{C_m}{M^{m-1}}=Q_0.
\end{\equ}
All inequalities are equalities. Therefore $h=M$ and $p=p_N$,
and $M^{(m-1)/2}\Phi$ is a ground-state solution.
The first bubble carries the entire limiting weighted mass, so
there is no remaining critical mass away from $p_N$.
The positive uniform lower bound for $H_j$ then shows that there is
no unweighted mass away from $p_N$, either. Since $H_j\to H_*$ uniformly
and $H_*(p_N)=M$, the total unweighted mass tends to
$Q_0/M=C_m/M^m$. This proves the measure convergence in
\eqref{eq:weak-coefficient-concentration}.
Small-energy regularity gives uniform convergence to zero on compact
sets outside $p_N$.

For completeness, the rescaled spinors in a fixed stereographic chart
are globally bounded in $L^q(\R^m)$, converge locally to $\Phi$, and
have total critical mass converging to that of $\Phi$.
The Brezis-Lieb lemma gives global strong $L^q$ convergence.
The rescaled coefficients are bounded and converge locally to $M$,
so the equations and dominated convergence give strong $L^{\hat q}$
convergence of their Dirac derivatives to $D\Phi$. The homogeneous
version of \eqref{eq:elliptic} upgrades this to strong
$\dot H^{1/2}(\R^m)$ convergence. These conclusions hold after
subsequence extraction for the profiles. Since every possible
extraction yields the same point $p_N$ and the same mass, the measure
convergence and divergence of the maximum norm hold for the whole
chosen sequence.
\end{proof}

\begin{remark}\label{rem:coefficient-regularity-counterexample}
The sequence has a common positive coefficient lower bound and total
critical mass tending to $C_m/M^m$. Hence it satisfies the hypotheses
of \cite[Theorem~1.1]{LT2026}. The uniform $C^2$ bound supplies the
local $C^1$ bounds required in \cite[Theorem~1.2]{LT2026}, which
further shows that the blow-up point is isolated simple. These results
provide additional local information about the sequence constructed
here.

The coefficients may also be chosen with a uniformly negative
Laplacian at every maximum. At $|y|=\rho_j$ the gradient vanishes,
and the coordinate Hessian is $4g''(1)\nu\ot\nu$, where
$\nu=y/\rho_j$. Uniform equivalence of the coordinate metric gives
$\De_{\ig_0}H_j\le-\eta<0$ on the maximum set for some fixed
$\eta>0$. The limiting coefficient satisfies
$\De_{\ig_0}H_*(p_N)=-ma<0$ as well.

The uniform $L^\infty$ estimate in \cite[Theorem~1.3]{LT2026}
holds near the first mass threshold for coefficients in a sufficiently
small $C^2$-neighborhood of a fixed coefficient satisfying this
Laplacian condition. Here $H_j$ does not converge to $H_*$ in $C^2$-topology.
We can still apply the theorem in \cite{LT2026} with respect to each $H_j$ as a fixed coefficient,
but the allowed mass increment and the resulting bound may then
depend on $j$. By \eqref{eq:weak-coefficient-concentration}, for any
fixed $\vr>0$, the unweighted mass of $\psi_j$ is at most
$C_m/M^m+\vr$ for all sufficiently large $j$, while
$\norm{\psi_j}_\infty\to\infty$. Thus the convergence in
\eqref{eq:weaker-coefficient-limits}, even together with a uniform
$C^2$-bound and a uniformly negative Laplacian at every global maximum,
does not exclude blow-up near the first mass threshold.
\end{remark}

\begin{remark}
For $m\ge3$, suppose that the infimum $E_{\mathrm{exc}}(m)$ of the
excited solution energies lies below $2E_m$.
By Proposition~\ref{prop:compactness-mod-conf}, a minimizing sequence
has a subsequence that converges strongly in $\msh$ to a nontrivial
solution after suitable conformal transformations. The limit has
energy $E_{\mathrm{exc}}(m)>E_m$ by Theorem~\ref{thm:gap}, so it is
an excited solution that attains the infimum. If $E_{\mathrm{exc}}(m)\ge2E_m$, then $\ka_m=1$ by
\eqref{eq:least-excited-definition2} and \eqref{eq:kappa-canonical}.
In this case the lower bound $d>2^{-1/(m-1)}H_{\max}$ in
\textup{($\mathrm{H}^*$)} already implies the pinching condition
\eqref{eq:pinching} in Theorem~\ref{thm:existence}.
\end{remark}

\medskip

\section*{AI Declaration Statement}
All mathematical ideas, arguments, and results presented in this manuscript are solely due
to the authors, and have been obtained without the use of generative AI. Chat GPT was partly exploited in the revision process, towards notational uniformization and streamlined exposition.

\end{document}